\documentclass[10pt]{article}

\usepackage{amssymb,amsmath,amsthm}

\theoremstyle{plain}
\newtheorem{theo}{Theorem}[section]
\newtheorem{lem}{Lemma}[section]
\newtheorem{co}{Corollary}[section]

\theoremstyle{definition}
\newtheorem{defn}{Definition}[section]

\theoremstyle{remark}
\newtheorem{rmk}{Remark}[section]
\newcommand{\psink}{\Psi_{{\bf n},k}}
\newcommand{\psinu}{\Psi_{{\bf n},1}}
\newcommand{\psind}{\Psi_{{\bf n},2}}

\newcommand{\tdt}{\tau_{2,3}}
\newcommand{\ttd}{\tau_{3,2}}
\newcommand{\sut}{s_{1,3}}

\newcommand{\intdu}{\int_{\Delta_1}}
\newcommand{\intdd}{\int_{\Delta_2}}
\newcommand{\intdt}{\int_{\Delta_3}}

\newcommand{\supp}[1]{{\rm supp\/}(#1)}

\begin{document}

\begin{center}
{\bf RATIO ASYMPTOTIC OF HERMITE-PAD\'{E} ORTHOGONAL POLYNOMIALS FOR NIKISHIN SYSTEMS. II}\\
{\bf Abey L\'opez Garc\'ia} \footnotemark \footnotetext{Dpto. de
Matem\'{a}ticas, Universidad Carlos III de Madrid, Avda.
Universidad 15, 28911 Legan\'{e}s, Madrid (SPAIN)
$<$ablopez@math.uc3m.es$>.$} \hspace{1cm} {\bf Guillermo L\'{o}pez
Lagomasino} \footnotemark \footnotetext{Dpto. de Matem\'{a}ticas,
Universidad Carlos III de Madrid, Avda. Universidad 15, 28911
Legan\'{e}s, Madrid (SPAIN) $<$lago@math.uc3m.es$>$}
\end{center}

\begin{abstract}
We prove ratio asymptotic for sequences of multiple orthogonal
polynomials with respect to a Nikishin system of measures
${\mathcal{N}}(\sigma_1,\ldots,\sigma_m)$ such that for each $k$,
the support of $\sigma_k$ consists of an interval
$\widetilde{\Delta}_k$, on which $\sigma_k^{\prime}
> 0$ almost everywhere, and a  set without accumulation points in
$\mathbb{R} \setminus \widetilde{\Delta}_k$.
\end{abstract}



\section{Introduction}
Let $s$ be a finite positive Borel measure supported on a bounded
interval $\Delta$ of the real line ${\mathbb{R}}$ such that
$s^{\prime} > 0$ almost everywhere on $\Delta$ and let $\{Q_n\},
n\in {\mathbb{Z}}_+,$ be the corresponding sequence of monic
orthogonal polynomials; that is, with leading coefficients equal
to one. In a series of two papers (see \cite{kn:Rak1} and
\cite{kn:Rak2}), E. A. Rakhmanov proved that under these
conditions
\begin{equation} \label{Rakteo}
\lim_{n \in {\mathbb{Z}}_+} \frac{Q_{n+1}(z)}{Q_n(z)} =
\frac{\varphi(z)}{\varphi^{\prime}(\infty)}, \qquad \mathcal{K}
\subset {\mathbb{C}} \setminus \Delta
\end{equation}
(uniformly on each compact subset of ${\mathbb{C}} \setminus
\Delta$), where $\varphi(z)$ denotes the conformal representation
of $\overline{\mathbb{C}} \setminus \Delta$ onto $\{w: |w| > 1\}$
such that $\varphi(\infty) = \infty$ and $\varphi^{\prime}(\infty)
>0$. This result attracted great attention because of its
theoretical interest within the general theory of orthogonal
polynomials and its applications to the theory of rational
approximation of analytic functions. Simplified proofs of
Rakhmanov's theorem may be found in \cite{kn:Rak3} and
\cite{kn:Nev1}.

This result has been extended in several directions. Orthogonal
polynomials with respect to varying measures (depending on the
degree of the polynomial) arise in the study of multipoint Pad\'e
approximation of Markov functions. In this context, in
\cite{kn:Gui1} and \cite{kn:Gui2}, an analogue of Rakhmanov's
theorem for such sequences of orthogonal polynomials was proved.
Recently, S. A. Denisov \cite{Denisov} (see also \cite{NT})
obtained a remarkable extension of Rakhmanov's result to the case
when the support of $s$ verifies $\supp s = \widetilde{\Delta}
\cup e \subset \mathbb{R}$, where $\widetilde{\Delta}$ is a
bounded interval, $e$ is a set without accumulation points  in
$\mathbb{R} \setminus \widetilde{\Delta} $, and $s^{\prime} > 0$
a.e. on $\widetilde{\Delta}$. A version for orthogonal polynomials
with respect to varying Denisov type measures was given in
\cite{BCL}.

Another direction of generalization is connected with multiple
orthogonal polynomials. These are polynomials whose orthogonality
relations are distributed between several measures. They appear as
the common denominator of Hermite-Pad\'e approximations of systems
of Markov functions. An interesting class of such systems is
formed by the so called Nikishin systems of functions introduced
in \cite{kn:Nikishin}. For Nikishin multiple orthogonal
polynomials a version of Rakhmanov's theorem was proved in
\cite{AptLopRoc}.

An elegant notation for Nikishin systems was proposed in
\cite{GonRakhSor}. Let $\sigma_1, \sigma_2$ be two finite Borel
measures with constant sign, whose supports $\supp{\sigma_1},$ $
\supp{\sigma_2}$ are contained in non intersecting intervals of
${\mathbb{R}}$. Set
\[ d \langle \sigma_1, \sigma_2 \rangle(x) = \int \frac{d
\sigma_2(t)}{x-t}\, d \sigma_1(x) = \widehat{\sigma}_2(x)d
\sigma_1(x)\,.
\]
This expression defines a new measure with constant sign whose
support coincides with that of $\sigma_1$. Whenever convenient, we
use the differential notation of a measure.

Let $\Sigma = (\sigma_1,\ldots,\sigma_m)$ be a system of finite
Borel measures on the real line with constant sign and compact
support  containing infinitely many points. Let $\mbox{Co}(\supp
{\sigma_k}) = \Delta_k$ denote the smallest interval which
contains $\supp {\sigma_k}$. Assume that
\[\Delta_k \cap \Delta_{k+1} =\emptyset\,,\qquad
k=1,\ldots,m-1\,.
\]
By definition, $S = (s_1,\ldots,s_m)=
{\mathcal{N}}(\sigma_1,\ldots,\sigma_m)$, where
\begin{equation} \label{eq:systema} s_1 = \sigma_1, \quad
s_2 = \langle \sigma_1,\sigma_2 \rangle, \ldots , \quad s_m =
\langle \sigma_1, \langle \sigma_2,\ldots,\sigma_m \rangle \rangle
\end{equation}
is called the {\em Nikishin system} of measures  generated by
$\Sigma$. The system $(\widehat{s}_1,\ldots,\widehat{s}_m)$ of
Cauchy transforms of a Nikishin system of measures gives a
Nikishin system of functions.

Fix a multi-index ${\bf n}=(n_1,\ldots,n_m)\in \mathbb{Z}_+^m$.
The polynomial $Q_{{\bf n}}(x)$ is called an ${\bf n}$-th {\em
multiple orthogonal polynomial} with respect to $S$ if it is not
identically equal to zero, $\deg {Q_{\bf n}} \leq |{\bf n}| = n_1
+ \cdots + n_m$, and
\begin{equation}            \label{eq:ortogonalidad}
\int Q_{{\bf n}}(x)x^\nu ds_k(x)=0, \quad \nu=0,\ldots,n_k-1,\quad
k=1,\ldots,m.
\end{equation}
In the sequel, we assume that $Q_{\bf n}$ is monic.

If (\ref{eq:ortogonalidad}) implies that $\deg Q_{{\bf n}} = |{\bf
n}|$, the multi--index ${\bf n}$ is said to be {\bf normal} and
the corresponding monic multiple orthogonal polynomial is uniquely
determined. In addition, if the zeros of $Q_{{\bf n}}$ are simple
and lie in the interior of $\mbox{Co}(\supp {\sigma_1}) $ the
multi--index is said to be {\bf strongly normal}. (In relation to
intervals of the real line the interior refers to the Euclidean
topology of ${\mathbb{R}}$.) For Nikishin systems with $m =1,2,3$,
all multi-indices are strongly normal (see \cite{FidLop}). An open
question is whether or not this is true for all $m \in
{\mathbb{N}}$. The best result when $m\geq 4$ is that all
\[ {\bf n} \in {\mathbb{Z}}_+^m(*) = \{{\bf n} \in {\mathbb{Z}}_+^m: \not
\exists \,\, 1 \leq i < j < k  \leq m\, ,\,\,\,\mbox{with}
\,\,\,n_i < n_j < n_k\}
\]
are strongly normal (see \cite{FidLop2}).

In \cite{AptLopRoc}, a Rakhmanov type theorem was proved for
Nikishin systems  such that $\sigma_k^{\prime} > 0$ a.e. on
$\mbox{Co}(\supp{\sigma_k}), k=1,\ldots,m,$ and sequences of
multi-indices contained in
\[ {\mathbb{Z}}_+^m(\circledast) = \{{\bf n} \in {\mathbb{Z}}_+^m:
1 \leq i < j \leq m \Rightarrow n_j \leq n_i +1\}\,.
\]
It is easy to see that ${\mathbb{Z}}_+^m(\circledast)\subset
{\mathbb{Z}}_+^m(*)$. Here, we assume that $\supp{\sigma_k} =
\widetilde{\Delta}_k \cup e_k, k=1,\ldots,m$, where
$\widetilde{\Delta}_k$ is a bounded interval of the real line,
$\sigma_k^{\prime} > 0$ a.e. on $\widetilde{\Delta}_k$, $e_k$ is a
set without accumulation points in $\mathbb{R} \setminus
\widetilde{\Delta}_k$, and the sequence of multi-indices on which
the limit is taken is in ${\mathbb{Z}}_+^m(*)$.

The proof of Theorem \ref{teo2} below uses the construction of so
called second type functions. This construction depends on the
relative value of the components of the multi-indices in
${\mathbb{Z}}_+^m(*)$ under consideration. A crucial step in our
study consists in proving an interlacing property for the zeros of
the second type functions corresponding to ``consecutive''
multi-indices (see Lemma \ref{entrelazamiento}). For this purpose,
we need to be sure that the second type functions are built using
the same procedure. To distinguish different classes of
multi-indices which respond for the same construction of second
type functions, we introduce the following definition.

\begin{defn}
Suppose that ${\bf n}=(n_1,\ldots,n_m)\in\mathbb{Z}_{+}^{m}$. Let
$\tau_{\bf n}$   denote the permutation of $ \{1,2,\ldots,m\}$
given by
\[\tau_{\bf n}(i)=j \quad \mbox{if }  \quad
\left\{
\begin{array}{cccc}
n_j > n_k & \mbox{for} & k<j, & k \not\in \{\tau_{\bf
n}(1),\ldots,\tau_{\bf n}(i-1)\} \\
n_j \geq n_k & \mbox{for} & k>j, & k \not\in \{\tau_{\bf
n}(1),\ldots,\tau_{\bf n}(i-1)\}
\end{array}
\right. \,.
\]
\end{defn}

In words, $\tau_{\bf n}(1)$ is the subindex of the first component
of ${\bf n}$ (from left to right) which is greater or equal than
the rest, $\tau_{\bf n}(2)$ is the subindex of the first component
which is second largest, and so forth. For example, if $n_1 \geq
\cdots \geq n_m$ then $\tau_{\bf n}$ is the identity.

Let $\tau$ denote a permutation of $\{1,2,\ldots,m\}$. Set
\[{\mathbb{Z}}_+^m(*,\tau) = \{{\bf n} \in {\mathbb{Z}}_+^m(*): \tau_{\bf n} =
\tau\}\,.
\]
Let ${\bf n} \in \mathbb{Z}_+^m$ and $l \in \{1,\ldots,m\}$.
Define
$$
{\bf n}_{l}:=(n_1,\ldots,n_{l-1},n_l+1,n_{l+1},\ldots,n_m)\,.
$$

Consider the $(m+1)$-sheeted Riemann surface
$$
\mathcal R=\overline{\bigcup_{k=0}^m \mathcal R_k} ,
$$
formed by the consecutively ``glued'' sheets
$$
\mathcal R_0:=\overline {\mathbb{C}} \setminus
\widetilde{\Delta}_1,\quad \mathcal R_k:=\overline {\mathbb{C}}
\setminus (\widetilde{\Delta}_k \cup
\widetilde{\Delta}_{k+1}),\,\, k=1,\dots,m-1,\quad \mathcal
R_m=\overline {\mathbb{C}} \setminus \widetilde{\Delta}_m,
$$
where the upper and lower banks of the slits of two neighboring
sheets are identified. Fix $l \in \{1,\ldots,m\}$. There exists a
conformal representation $G^{(l)}$ of $\mathcal{R}$ onto
$\overline{\mathbb{C}}$ such that
\[G^{(l)}(z) = z + \mathcal{O}(1)\,,\,\, z \to \infty^{(0)}, \quad
G^{(l)}(z) = C/z + \mathcal{O}(1/z^2)\,,\,\, z \to \infty^{(l)}.
\]
By $G^{(l)}_k$ we denote the branch of $G^{(l)}$ on
$\mathcal{R}_k$.

\begin{theo} \label{teo2}
Let $S = {\mathcal{N}}(\sigma_1,\ldots,\sigma_m)$ be a Nikishin
system with $ \supp{\sigma_k}  = \widetilde{\Delta}_k \cup e_k,
k=1,\ldots,m$, where $\widetilde{\Delta}_k$ is a bounded interval
of the real line, $\sigma_k^{\prime} > 0$ a.e. on
$\widetilde{\Delta}_k$, and $e_k$ is a set without accumulation
points in $\mathbb{R} \setminus \widetilde{\Delta}_k$. Let
$\Lambda \subset {\mathbb{Z}}_+^m(*)$ be an infinite sequence of
distinct multi-indices. Let us assume that there exists $l\in
\{1,\ldots,m\}$ and a fixed permutation $\tau$ of $\{1,\ldots,m\}$
such that for all ${\bf n} \in \Lambda$ we have that ${\bf n},{\bf
n}_l \in \mathbb{Z}_{+}^{m}(*,\tau)$  and
$\displaystyle{\max_{{\bf n} \in \Lambda }(\max_{k
=1,\ldots,m}mn_{k} -  |{\bf n}| ) }< \infty.$ Then,
\begin{equation} \label{eq:xe} \lim_{{\bf n}\in
{\Lambda}}\frac{Q_{{\bf n}_{l} }(z)}{Q_{{\bf n} }(z)}=
G_0^{(l)}(z), \qquad  \mathcal{K}  \subset {\mathbb{C}} \setminus
\supp{\sigma_1} \,.
\end{equation}
\end{theo}

When $m=1$ this result reduces to Denisov's version of Rakhmanov's
theorem. The proof of Theorem \ref{teo2} follows the guidelines
employed in \cite{AptLopRoc} but it is technically more
complicated because of the more general assumptions on the
measures and the sequence of multi-indices.

Let ${\bf 1} = (1,\ldots,1).$ An immediate consequence of Theorem
\ref{teo2} is
\begin{co} \label{cor1} Let $S = {\mathcal{N}}(\sigma_1,\ldots,\sigma_m)$ be a Nikishin
system with $ \supp{\sigma_k}  = \widetilde{\Delta}_k \cup e_k,
k=1,\ldots,m$, where $\widetilde{\Delta}_k$ is a bounded interval
of the real line, $\sigma_k^{\prime} > 0$ a.e. on
$\widetilde{\Delta}_k$, and $e_k$ is a set without accumulation
points in $\mathbb{R} \setminus \widetilde{\Delta}_k$. Let
$\Lambda \subset {\mathbb{Z}}_+^m(*)$ be an infinite sequence of
distinct multi-indices such that $\displaystyle{\max_{{\bf n} \in
\Lambda }(\max_{k =1,\ldots,m}mn_{k} -  |{\bf n}| ) }< \infty.$
Then,
\begin{equation} \label{eq:xe1} \lim_{{\bf n}\in
{\Lambda}}\frac{Q_{{\bf n +1}}(z)}{Q_{{\bf n} }(z)}= \prod_{l=1}^m
G_0^{(l)}(z), \qquad  \mathcal{K}  \subset {\mathbb{C}} \setminus
\supp{\sigma_1} \,.
\end{equation}
\end{co}

The paper is organized as follows. In Section 2 we introduce and
study an auxiliary system of second type functions. An interlacing
property for the zeros of the polynomials $Q_{\bf n}$ and of the
second type functions is proved in Section 3. Using the
interlacing property of zeros and results on ratio and relative
asymptotic of polynomials orthogonal with respect to varying
measures, in Section 4 a system of boundary value problems is
derived which implies the existence of limit in (\ref{eq:xe}).
Actually, a more general result is proved which also contains the
ratio asymptotic of the second type functions.

\section{{\large Functions of second type and orthogonality
properties}}

Fix ${\bf n} =(n_1,\ldots,n_m)\in \mathbb{Z}_{+}^{m}(\ast)$ and
consider $Q_{\bf n}$ the ${\bf n}$-th multi-orthogonal polynomial
with respect to a Nikishin system $S=\mathcal{N}(\Sigma),\,
\Sigma=(\sigma_1,\ldots,\sigma_m).$ For short, in the sequel we
denote $\Delta_k = \mbox{Co}(\supp{\sigma_k}), k=1,\ldots,m.$
Inductively, we define functions of second type
$\psink,\,k=0,1,\ldots,m$, systems of measures $\Sigma^{k} =
(\sigma_{k+1}^k,\ldots,\sigma_m^k), k = 0,1,\ldots,m-1,
\mbox{Co}(\supp{\sigma_j^k}) \subset \Delta_j,$ which generate
Nikishin systems, and multi-indices ${\bf n}^k \in
\mathbb{Z}_+^{m-k}(*)  , k=0,\ldots,m-1$. Take $\Psi_{{\bf
n},0}=Q_{\bf n}, {\bf n}^0={\bf n},$ and $\Sigma^0=\Sigma$.

Suppose that ${\bf n}^k=(n_{k+1}^{k},\ldots,n_{m}^{k})$,
$\Sigma^{k}=(\sigma_{k+1}^{k},\ldots,\sigma_{m}^{k})$ and $\psink$
have already been defined, where $0\leq k\leq m-2$. Let
\[{\bf n}^{k+1}=(n_{k+2}^{k+1},\ldots,n_{m}^{k+1})\in\mathbb{Z}_{+}^{m-k-1}(\ast)\]
be the multi-index obtained  deleting from ${\bf n}^k$ the first
component $n_{r_k}^k$ which verifies
\[n_{r_k}^k=\max\{n_j^k:\, k+1\leq j\leq m\}.\]
The components of ${\bf n}^{k +1}$ and ${\bf n}^{k}$ are related
as follows:
\[n_{k+1}^{k}=n_{k+2}^{k+1},\ldots,n_{r_k-1}^k=n_{r_k}^{k+1},n_{r_k+1}^k=n_{r_k+1}^{k+1},\ldots,n_{m}^{k}=n_{m}^{k+1}.\]
Denote
\begin{equation}\label{eq1}
\Psi_{{\bf
n},k+1}(z)=\int_{\Delta_{k+1}}\frac{\psink(x)}{z-x}\,ds_{r_k}^k(x)\,,
\end{equation}
where $s_{r_k}^k=\langle
\sigma_{k+1}^{k},\ldots,\sigma_{r_k}^{k}\rangle$ is the
corresponding component of the Nikishin system
$S^k=\mathcal{N}(\Sigma^k)=(s_{k+1}^k,\ldots,s_m^k)$.

In order to define $\Sigma^{k+1}$ we introduce the following
notation. Set
\[s_{i,j}^k=\langle
\sigma_{i}^k,\ldots,\sigma_j^k\rangle,\qquad k+1\leq i\leq j\leq
m,\] where  $\sigma^k_i \in \Sigma^{k}$. In page 390 of
\cite{Krein} it is proved that there exists a finite measure
$\tau_{i,j}^k$ with constant sign such that
\[\mbox{Co}(\supp {\tau_{i,j}^k})\subset \mbox{Co}(\supp{s_{i,j}^k}) \]
\[\frac{1}{\widehat{s}_{i,j}^k(z)}=l_{i,j}^k(z)+\widehat{\tau}_{i,j}^k(z)\]
where $l_{i,j}^k$ is a certain polynomial of degree 1.  That
$\mbox{Co}(\supp {s_{i,j}^k}) \subset \Delta_i$ easily follows by
induction. We wish to remark that the continuous part of
$\supp{s^k_{i,j}}$ and $\supp{\tau^k_{i,j}}$ coincide, but not
their isolated parts. In fact, zeros of $\widehat {s}^k_{i,j}$ on
$\Delta_i$ (there is one such zero between two consecutive mass
points of ${s}^k_{i,j}$) become poles of $\widehat{\tau}^k_{i,j}$
(mass points of ${\tau}^k_{i,j}$).

Suppose that $r_k=k+1$. In this case, we take
\[\Sigma^{k+1}=(\sigma_{k+2}^k,\ldots,\sigma_{m}^k) =
(\sigma_{k+2}^{k+1},\ldots,\sigma_{m}^{k+1})\] deleting the first
measure of $\Sigma^k$. If $r_k\geq k+2,$ then $\Sigma^{k+1}$ is
defined by
\[(\tau_{k+2,r_k}^k,\widehat{s}_{k+2,r_k}^k
d\tau_{k+3,r_k}^k,\ldots,\widehat{s}_{r_k-1,r_k}^k
d\tau_{r_k,r_k}^k,\widehat{s}_{r_k,r_k}^k
d\sigma_{r_k+1}^k,\sigma_{r_k+2}^k,\ldots,\sigma_m^k)\,,
\]
where $\mbox{Co}(\supp{\sigma^{k+1}_j}) \subset \Delta_j,
j=k+2,\ldots,m.$ Any two consecutive measures in the system
$\Sigma^{k+1}$ are supported on disjoint intervals; therefore, $
\Sigma^{k+1}$ generates a Nikishin system. To conclude we define
\[
\Psi_{{\bf n},m}(z)=\int_{\Delta_{m}}\frac{\Psi_{{\bf
n},m-1}(x)}{z-x}\,ds_{m}^{m-1}(x)\,.
\]

If   $n_1 \geq \cdots \geq n_m$, we have that ${\bf n}^k =
(n_{k+1},\ldots,n_m), \Sigma^k = (\sigma_{k+1},\ldots,\sigma_m)$
and $\Psi_{{\bf n},k}(z) = \int_{\Delta_k} \frac{\Psi_{{\bf
n},k-1}(x)}{z-x} d\sigma_k(x), k=1,\ldots,m$. Basically, this is
the situation considered in \cite{AptLopRoc}.

To fix ideas let us turn our attention to the cases $m=2$ and
$m=3$. We denote by $\mathcal{C}(f;\mu)$ the Cauchy transform of
$fd\mu$; that is,
\[\mathcal{C}(f;\mu)(z)=\int\frac{f(x)}{z-x}\,d\mu(x)\,.\] In
the following tables, we omit the line corresponding to $k=0$
because by definition $\Sigma^0=\Sigma$, $\Psi_{{\bf n},0}=Q_{\bf
n}$ and ${\bf n}^0={\bf n}$.

\begin{table}[h] \caption{m=2} \label{tabla1}
\begin{center} \footnotesize
\begin{tabular}{|c|c|c|c|c|c|}
\hline
 ${m =2}$ &   ${k}$ &  ${r_{k-1}}$ &
 ${\Psi_{{\bf n},k}}$ &
$ {\Sigma^k} $&    ${{\bf n}^k} $ \\ \hline $n_1 \geq n_2 $ & $1$
& $1$ & ${\mathcal{C}}(Q_{{\bf n}};\sigma_1)$ & $(\sigma_2)$ &
$(n_2)$ \\ \hline $ n_1 < n_2$ & $1$ & $2$ &
${\mathcal{C}}(Q_{{\bf
n}};\langle \sigma_1, \sigma_2 \rangle)$ & $(\tau_2)$ &   $(n_1)$ \\
\hline
\end{tabular}
\end {center}
\end{table}

\begin{table}[h] \caption{$m=3$} \label{tabla2}
\begin{center} \footnotesize
\begin{tabular}{|c|c|c|c|c|c|}
\hline ${m =3}$ & $  {k}$ & $ {r_{k-1}}$ &
 ${\Psi_{{\bf n},k}} $ &
 ${\Sigma^k}$ &   $ {{\bf n}^k} $ \\ \hline
$n_1 \geq n_2 \geq n_3 $& $1$ & $1  $& ${\mathcal{C}}(Q_{{\bf
n}};\sigma_1)$ & $ (\sigma_2,\sigma_3)$ &   $(n_2,n_3)$ \\
\hline $\mbox{} $&$ 2 $&$ 2 $&$ {\mathcal{C}}(\Psi_{{\bf
n},1};\sigma_2)$ & $( \sigma_3)$ & $  ( n_3) $ \\ \hline $n_1 \geq
n_3
> n_2$ & $1$ & $1$ &${\mathcal{C}}(Q_{{\bf n}};\sigma_1)$ &$
(\sigma_2,\sigma_3)$ &$   (n_2,n_3) $\\
\hline $\mbox{} $&$ 2 $&$ 3 $&$ {\mathcal{C}}(\Psi_{{\bf
n},1};\langle\sigma_2,\sigma_3\rangle) $&$ ( \tau_3)$ &$   ( n_2)
$\\ \hline $n_2 > n_1 \geq n_3 $&$ 1 $&$ 2 $&$
{\mathcal{C}}(Q_{{\bf n}};\langle\sigma_1,\sigma_2\rangle) $&$
(\tau_2,\langle\sigma_3,\sigma_2\rangle) $&$  (n_1,n_3)  $\\
\hline $\mbox{} $ & $2$ & $2$ &$ {\mathcal{C}}(\Psi_{{\bf n},1};
\tau_2), $&$ (\langle \sigma_3,\sigma_2 \rangle) $&$  ( n_3)  $\\
\hline $ n_2 \geq n_3 > n_1 $& $1$ &$ 2$ &$ {\mathcal{C}}(Q_{{\bf
n}};\langle\sigma_1,\sigma_2\rangle) $&$
(\tau_2,\langle\sigma_3,\sigma_2\rangle) $&$   (n_1 ,n_3) $ \\
\hline $\mbox{} $&$ 2 $&$ 3 $&$ {\mathcal{C}}(\Psi_{{\bf n},1};
\langle \tau_2, \sigma_3,\sigma_2\rangle) $& $(\tau_{3,2}) $&$   (
n_1)  $\\ \hline $n_3 > n_1 \geq n_2 $&$ 1$ &$ 3$ &$
{\mathcal{C}}(Q_{{\bf
n}};\langle\sigma_1,\sigma_2,\sigma_3\rangle)$ &$
(\tau_{2,3},\langle\tau_3,\sigma_2,\sigma_3\rangle) $&$  (n_1,n_2)
$\\ \hline $\mbox{} $& $2$ &$ 2$ &$ {\mathcal{C}}(\Psi_{{\bf
n},1}; \tau_{2,3})$ & $(\langle\tau_3,\sigma_2,\sigma_3\rangle)$ &
  $( n_2)$ \\ \hline
\end{tabular}
\end {center}
\end{table}

In Theorem 2 of \cite{FidLop2} it was proved that the functions
$\psink$ verify the following orthogonality relations. For each
$k=0,1,\ldots,m-1,$
\begin{equation}\label{eq2}
\int_{\Delta_{k+1}}x^{\nu}\psink(x)\,ds_{i}^k(x)=0,\quad
\nu=0,1,\ldots,n_{i}^k-1,\quad i=k+1,\ldots,m\,,
\end{equation}
where $s_i^k=\langle \sigma_{k+1}^k,\ldots,\sigma_i^k\rangle$.

We wish to underline that since
$\mathbb{Z}_{+}^{2}(\ast)=\mathbb{Z}_{+}^{2}$, all multi-indices
with two components have associated functions of second type.
However, for $m=3$ the case $n_1<n_2<n_3$ has not been considered
(see Table \ref{tabla2}). The rest of this section will be devoted
to the construction of certain functions $\psink$ for this case
and to the proof of the orthogonality relations they satisfy. We
use the following auxiliary result.
\begin{lem}\label{lemaint} Let $s_{3,2} = \langle
\sigma_3,\sigma_2\rangle$. Then
\begin{equation}\label{eq4}
\int_{\Delta_2}\frac{\widehat{s}_{3,2}(x)}{\widehat{\sigma}_{3}(x)}\,
\frac{d\tau_{2,3}(x)}{(z-x)} + C_1
=\frac{\widehat{\sigma}_2(z)}{\widehat{s}_{2,3}(z)} \,,\quad z \in
\mathbb{C}\setminus \supp{\sigma_2} \,,
\end{equation}
where $C_1 = \sigma_2(\Delta_2)/s_{2,3}(\Delta_2)$.
\end{lem}

\begin{proof} We employ two useful relations. The first one is
\begin{equation}\label{eq6}
\widehat{\sigma}_2(\zeta)\,\widehat{\sigma}_3(\zeta)=
\widehat{s}_{2,3}(\zeta)+\widehat{s}_{3,2}(\zeta),\quad
\zeta\in\mathbb{C}\setminus(\supp{\sigma_2}\cup\supp{\sigma_3})\,.
\end{equation}
The proof is straightforward and may be found in Lemma 4 of
\cite{FidLop}. The second one was mentioned above and states that
there exists a polynomial $l_{2,3}$ of degree $1$ and a measure
$\tau_{2,3}$ such that
\begin{equation}\label{eq5}
\frac{1}{\widehat{s}_{2,3}(z)}=\widehat{\tau}_{2,3}(z)+l_{2,3}(z),\quad
z\in \mathbb{C}\setminus \supp{\sigma_{2}}\,.
\end{equation}

Notice that
\[
\frac{\widehat{\sigma}_2(z)}{\widehat{s}_{2,3}(z)} - C_1 =
\mathcal{O}\left( \frac{1}{z}\right) \in
\mathcal{H}(\overline{\mathbb{C}}\setminus \Delta_2)
\]
Let $\Gamma$ be a positively oriented smooth closed Jordan curve
such that $\Delta_2$ and $\{z\} \cup \Delta_3$ lie on the bounded
and unbounded connected components, respectively, of
$\mathbb{C}\setminus\Gamma$.  By Cauchy's integral formula, we
have
\[ \frac{\widehat{\sigma}_2(z)}{\widehat{s}_{2,3}(z)}-
C_1 = \frac{1}{2\pi i} \int_{\Gamma}
\left(\frac{\widehat{\sigma}_2(\zeta)}{\widehat{s}_{2,3}(\zeta)}-
C_1\right) \frac{d\zeta}{z - \zeta} = \frac{1}{2\pi i}
\int_{\Gamma}
 \frac{\widehat{\sigma}_2(\zeta)}{\widehat{s}_{2,3}(\zeta)}
 \frac{d\zeta}{z - \zeta}\,.
\]
Multiply and divide the expression under the last integral sign by
$\widehat{\sigma}_3$ and use (\ref{eq6}) to obtain
\[\frac{\widehat{\sigma}_2(z)}{\widehat{s}_{2,3}(z)}-
C_1 = \frac{1}{2\pi i}\int_{\Gamma} \frac{\widehat{s}_{2,3}(\zeta)
+
\widehat{s}_{3,2}(\zeta)}{\widehat{\sigma}_3(\zeta)\widehat{s}_{2,3}(\zeta)}
 \frac{d\zeta}{z - \zeta} = \frac{1}{2\pi i}\int_{\Gamma}
\frac{
\widehat{s}_{3,2}(\zeta)}{\widehat{\sigma}_3(\zeta)\widehat{s}_{2,3}(\zeta)}
 \frac{d\zeta}{z - \zeta}\,.
\]
Taking account of (\ref{eq5}) it follows that
\[\frac{\widehat{\sigma}_2(z)}{\widehat{s}_{{2,3}}(z)}-
C_1 =\frac{1}{2\pi i}\int_{\Gamma} \frac{
\widehat{s}_{{3,2}}(\zeta)}{\widehat{\sigma}_3(\zeta)}\frac{(\widehat{\tau}_{2,3}(\zeta)
+ l_{2,3}(\zeta))d\zeta}{z - \zeta} =\frac{1}{2\pi i}\int_{\Gamma}
\frac{ \widehat{s}_{{3,2}}(\zeta)}{\widehat{\sigma}_3(\zeta)}
  \frac{\widehat{\tau}_{2,3}(\zeta)d\zeta}{z - \zeta}\,.
\]
Now, substitute $\widehat{\tau}_{2,3}(\zeta)$ by its integral
expression and use the Fubini  and Cauchy  theorems to obtain
\[\frac{\widehat{\sigma}_2(z)}{\widehat{s}_{{2,3}}(z)}-
C_1 = \int \frac{1}{2\pi i} \int_{\Gamma} \frac{
\widehat{s}_{{3,2}}(\zeta)}{\widehat{\sigma}_3(\zeta)(z - \zeta)}
  \frac{d\zeta}{\zeta - x} d\tau_{2,3}(x) = \int
\frac{ \widehat{s}_{{3,2}}(x)}{\widehat{\sigma}_3(x)}
\frac{d\tau_{2,3}(x)}{z - x}\,,
\]
which is what we set out to prove.
\end{proof}

We are ready to define  the functions of second type and to prove
the orthogonality properties they verify for multi-indices with 3
components not in $\mathbb{Z}_{+}^3(\ast)$ (with $n_1 < n_2 <
n_3$).

\begin{lem}\label{teoremaeq}
Fix ${\bf n}=(n_1,n_2,n_3)\in\mathbb{Z}_{+}^3$ where $n_1<n_2<n_3$
and consider $Q_{\bf n}$ the {\bf n}-th orthogonal polynomial
associated to a Nikishin system
$S=(s_1,s_2,s_3)=\mathcal{N}(\sigma_1,\sigma_2,\sigma_3)$ . Set
$\Psi_{{\bf n},0}=Q_{\bf n}$,
\begin{equation}\label{eq8}
\Psi_{{\bf n},1}(z)=\int_{\Delta_1}\frac{Q_{\bf
n}(x)}{z-x}\,d\,s_{1,3}(x)\,,
\end{equation}
\begin{equation}\label{eq9}
\Psi_{{\bf n},2}(z)=\int_{\Delta_2}\frac{\Psi_{{\bf
n},1}(x)}{z-x}\,\frac{\widehat{s}_{3,2}(x)}{\widehat{\sigma}_3(x)}\,d\,\tdt(x)\,.
\end{equation}
Then
\begin{equation}\label{eq10}
\int_{\Delta_1}t^{\nu}\,\Psi_{{\bf
n},0}(t)\,d\,s_{1,j}(t)=0,\qquad 0\leq\nu\leq n_j-1,\quad 1\leq
j\leq 3
\end{equation}
\begin{equation}\label{eq11}
\int_{\Delta_2}t^{\nu}\,\Psi_{{\bf n},1}(t)\,d\,\tdt(t)=0,\qquad
0\leq\nu\leq n_1-1
\end{equation}
\begin{equation}\label{eq12}
\int_{\Delta_2}t^{\nu}\,\Psi_{{\bf
n},1}(t)\,\frac{\widehat{s}_{3,2}(t)}{\widehat{\sigma}_3(t)}\,d\,\tdt(t)=0,\qquad
0\leq\nu\leq n_2-1
\end{equation}
\begin{equation}\label{eq13}
\int_{\Delta_3}t^{\nu}\,\Psi_{{\bf
n},2}(t)\,\frac{\widehat{s}_{2,3}(t)}{\widehat{\sigma}_2(t)}\,d\,\ttd(t)=0,\qquad
0\leq\nu\leq n_1-1.
\end{equation}
\end{lem}

\begin{rmk}
The measure $ {\widehat{s}_{3,2}}d\tau_{2,3}/{\widehat{\sigma}_3}
$ supported on $\Delta_2$ cannot be written in the form
$\langle\tdt,\mu\rangle$ for some measure $\mu$ supported on
$\Delta_3$, so there is no $\Sigma^1$ and $S^1$ in this case.
\end{rmk}

\begin{proof} The relations (\ref{eq10}) follow directly from the definition of
$Q_{\bf n}$. Let us justify (\ref{eq11}) and (\ref{eq12}).

For $0 \leq \nu \leq n_1-1 (\leq n_3 - 3)$, applying Fubini's
theorem,
\[\intdd t^{\nu}\,\Psi_{{\bf n},1}(t)\,d\tdt(t)=\intdd t^{\nu} \intdu \frac{Q_{\bf n}(x)}{t-x}\,d\sut(x)\,d\tdt(t)\]
\[=\intdu Q_{\bf n}(x) \intdd
\frac{t^{\nu}-x^{\nu}+x^{\nu}}{t-x}\,d\,\tdt(t)\,d\,\sut(x)\]
\[=\intdu Q_{\bf n}(x)\,p_{\nu}(x)\,d\,s_{1,3}(x)-\intdu
x^{\nu}Q_{\bf n}(x)\,\widehat{\tau}_{2,3}(x)\,d\,\sut(x)\,,\]
where $p_{\nu}(x)=\intdd \frac{t^{\nu}-x^{\nu}}{t-x}\,d\tdt(t)$ is
a polynomial of degree at most $n_1-2$. Since
$d\sut(x)=\widehat{s}_{2,3}(x)d\sigma_1(x)$ and
$\widehat{\tau}_{2,3}(x)\,\widehat{s}_{2,3}(x)=1-l_{2,3}(x)\,\widehat{s}_{2,3}(x)$,
the measure $\widehat{\tau}_{2,3}(x)\,d\sut(x)$ is equal to
$d\,\sigma_1(x)-l_{2,3}(x)\,d\sut(x)$. Therefore, applying
(\ref{eq10}) both integrals vanish and we obtain (\ref{eq11}).
Actually, we only needed that $n_1 \leq n_3 -1$.

If $0\leq\nu\leq n_2-1 (\leq n_3-2)$,
\[\intdd t^{\nu}\,\psinu
(t)\,\frac{\widehat{s}_{3,2}(t)}{\widehat{\sigma}_{3}(t)}\,d\tdt(t)=\intdd
t^{\nu}\frac{\widehat{s}_{3,2}(t)}{\widehat{\sigma}_{3}(t)}\intdu\frac{Q_{\bf
n}(x)}{t-x}\,d\sut(x)\,d\tdt(t)\]
\[=\intdu Q_{\bf n}(x)\intdd
\frac{t^{\nu}-x^{\nu}+x^{\nu}}{t-x}\,\frac{\widehat{s}_{3,2}(t)}{\widehat{\sigma}_{3}(t)}\,d\tdt(t)\,d\sut(x)\]
\[=\intdu Q_{\bf n}(x)\,x^{\nu}\intdd
\frac{\widehat{s}_{3,2}(t)}{\widehat{\sigma}_{3}(t)}\,\frac{d\tdt(t)}{t-x}d\sut(x)\]
By Lemma ~\ref{lemaint}, the last expression is equal to
\[C_1\,\intdu Q_{\bf n}(x)\,x^{\nu}\,d\sut(x)-\intdu Q_{\bf n}(x)\,x^{\nu}\,
\frac{\widehat{\sigma}_{2}(x)}{\widehat{s}_{2,3}(x)}\,d\sut(x)\]
\[=-\intdu Q_{\bf n}(x)\,x^{\nu}\,d s_{1,2}(x)=0\]
taking into account that
$d\sut(x)=\widehat{s}_{2,3}(x)\,d\sigma_{1}(x)$ and (\ref{eq10}).
This proves (\ref{eq12}). It would have been sufficient to require
$n_2 \leq n_3.$

Let us prove (\ref{eq13}). Take $0\leq \nu\leq n_1-1$, we have
\[\intdt t^{\nu}\,\psind (t)\,\frac{\widehat{s}_{2,3}(t)}{\widehat{\sigma}_{2}(t)}d\ttd(t)=
\intdt t^{\nu}\intdd
\frac{\psinu(x)}{t-x}\frac{\widehat{s}_{3,2}(x)}{\widehat{\sigma}_{3}(x)}d\tdt(x)
\frac{\widehat{s}_{2,3}(t)}{\widehat{\sigma}_{2}(t)}d\ttd(t)\]
\[=\intdd \psinu (x)\,\frac{\widehat{s}_{3,2}(x)}{\widehat{\sigma}_{3}(x)}\,\intdt
\frac{t^{\nu}-x^{\nu}+x^{\nu}}{t-x}\,\frac{\widehat{s}_{2,3}(t)}{\widehat{\sigma}_{2}(t)}d
\tau_{3,2}(t)\,d\tau_{2,3}(x) \]
\[= \intdd \,p_{\nu}(x)\,\psinu
(x)\,\frac{\widehat{s}_{3,2}(x)}{\widehat{\sigma}_{3}(x)}\,d\tau_{2,3}(x)\]
\[+\intdd\frac{\psinu(x)\,x^{\nu}\,\widehat{s}_{3,2}(x)}{\widehat{\sigma}_{3}(x)}\,\intdt
\frac{\widehat{s}_{2,3}(t)}{\widehat{\sigma}_{2}(t)}\frac{d\ttd(t)}{t-x}d\tdt(x)\]
where $p_{\nu}(x)$ is the polynomial defined by
\[ \intdt
\frac{t^{\nu}-x^{\nu}}{t-x}\,\frac{\widehat{s}_{2,3}(t)}{\widehat{\sigma}_{2}(t)}\,d\ttd
(t),\]  of degree   $\leq n_1-2$. Applying (\ref{eq12}), the first
integral after the last equality equals zero since $n_1 < n_2$
(though $n_1 \leq n_2 +1$ would have been sufficient). If we
interchange the sub-indices $2$ and $3$ in Lemma ~\ref{lemaint},
we obtain
\begin{equation}\label{eq14}
\intdt
\frac{\widehat{s}_{2,3}(t)}{\widehat{\sigma}_{2}(t)}\,\frac{d\ttd(t)}{t-x}=
-\frac{\widehat{\sigma}_{3}(x)}{\widehat{s}_{3,2}(t)}+ C_2 \,,
\end{equation}
where $C_2 =  \sigma_{3}(\Delta_3)/s_{3,2}(\Delta_3).$ Therefore,
using (\ref{eq14}), (\ref{eq12}) and (\ref{eq11}), it follows that
\[\intdd\frac{\psinu(x)\,x^{\nu}\,\widehat{s}_{3,2}(t)}{\widehat{\sigma}_{3}(x)}
\intdt\frac{\widehat{s}_{2,3}(t)}{\widehat{\sigma}_{2}(t)}\frac{d\ttd(t)}{t-x}\,d\tdt(x)
\]
\[=\intdd
\psinu(x)\,x^{\nu}\,\frac{\widehat{s}_{3,2}(t)}{\widehat{\sigma}_{3}(x)}
\left(C_2 -
\frac{\widehat{\sigma}_{3}(x)}{\widehat{s}_{3,2}(t)}\right)\,d\tdt(x)=0\,,\]
since $n_1 \leq n_2.$ This completes the proof.
\end{proof}

\section{{\large Interlacing property of zeros and varying measures}}

As we have pointed out,   from  the definition
${\mathbb{Z}}_+^m(*) = {\mathbb{Z}}_+^m, m=1,2$. We have
introduced adequate functions of second type also when $m=3$ and
$n_1 < n_2 <n_3$ which were the only multi-indices initially not
in ${\mathbb{Z}}_+^3(*)$. To unify notation, in the rest of the
paper we will consider that ${\mathbb{Z}}_+^3(*) =
{\mathbb{Z}}_+^3$.

In this section, we show that for ${\bf n} \in \mathbb{Z}_+^m(*),
m \in \mathbb{N}$, the functions $\Psi_{{\bf n},k},
k=0,\ldots,m-1,$ have exactly $|{\bf n}^k|$ simple zeros in the
interior of $\Delta_{k+1}$ and no other zeros on $\mathbb{C}
\setminus \Delta_{k}$. The zeros of ``consecutive'' $\Psi_{{\bf
n},k}$ satisfy an interlacing property. These properties are
proved in Lemma \ref{entrelazamiento} below which complements
Theorem 2.1 (see also Lemma 2.1) in \cite{AptLopRoc} and
substantially enlarges the class of multi-indices for which it is
applicable. The concept of AT system is crucial in its proof.

\begin{defn}
Let $(\omega_1,\omega_2,\ldots,\omega_m)$ be a collection of
functions which are analytic on a neighborhood of an interval
$\Delta$. We say that it forms an AT-system for the multi-index
${\bf n}=(n_1,n_2,\ldots,n_m)$ on $\Delta$ if whenever one chooses
polynomials $P_{n_1},\ldots,P_{n_m}$ with $\deg(P_{n_j})\leq
n_j-1$, not all identically equal to zero, the function
\[P_{n_1}(x)\,w_1(x)+\cdots+P_{n_m}(x)\,w_m(x)\]
has at most $|{\bf n}|-1$ zeros on $\Delta$, counting
multiplicities. $(\omega_1,\ldots,\omega_m)$ is an AT-system on
$\Delta$ if it is an AT-system on that interval for all ${\bf n}
\in \mathbb{Z}_{+}^m$.
\end{defn}

Theorem 1 of \cite{FidLop} (for $m=3$) and Theorem 1 of
\cite{FidLop2} prove the following.

\begin{lem} \label{lm:AT}
Let $(s_1,\ldots,s_{m-1}) =
\mathcal{N}(\sigma_1,\ldots,\sigma_{m-1}), m \geq 2,$ be a
Nikishin system of $m-1$ measures. Then
$(1,\widehat{s}_{1},\ldots,\widehat{s}_{{m-1}})$ forms an AT
system on any interval $\Delta$ disjoint from $\Delta_1$ with
respect to any ${\bf n} \in \mathbb{Z}_+^m(*)$.
\end{lem}

Recall that ${\bf n}_l$ denotes the multi-index obtained adding
$1$ to the $l$th component of ${\bf n}$.

\begin{lem} \label{entrelazamiento} Let $S =
{\mathcal{N}}(\sigma_1,\ldots,\sigma_m)$ be a Nikishin system. Let
${\bf n} \in \mathbb{Z}_+^m(*), m \in \mathbb{N}$, then for each
$k=0,\ldots,m-1,$ the function $\Psi_{{\bf n},k}$ has exactly
$|{\bf n}^k|$ simple zeros in the interior of $\Delta_{k+1}$ and
no other zeros on $\mathbb{C} \setminus \Delta_{k}$. Let $I$
denote the closure of any one of the connected components of
$\Delta_{k+1} \setminus \supp{\sigma_{k+1}^{k}},$ then $\Psi_{{\bf
n},k}$ has at most one simple zero on $I$. Assume that
$l\in\{1,2,\ldots,m\}$ is such that ${\bf n}, {\bf
n}_l\in\mathbb{Z}_{+}^{m}(\ast,\tau)$ for a fixed permutation
$\tau$. Then, for each $k \in \{0,\ldots,m-1\}$ between two
consecutive zeros of $\Psi_{{\bf n}_l,k}$ lies exactly one zero of
$\Psi_{{\bf n},k}$ and viceversa (that is, the zeros of
$\Psi_{{\bf n}_l,k}$ and $\Psi_{{\bf n},k}$ on $\Delta_{k+1}$
interlace).
\end{lem}
\begin{proof}   Assume that ${\bf n},
{\bf n}_l\in\mathbb{Z}_{+}^{m}(\ast,\tau)$. We claim that for any
real constants $A, B, |A| +|B| >0,$ and $k \in
\{0,1,\ldots,m-1\}$, the function
\[G_{{\bf n},k}(x) = A\Psi_{{\bf n},k}(x) +
B\Psi_{{\bf n}_l,k}(x)
\]
has at most $|{\bf n}^k| + 1  $ zeros in ${\mathbb{C}} \setminus
\Delta_k$ (counting multiplicities) and at least $|{\bf n}^k|$
simple zeros in the interior of $\Delta_{k+1} \,(\Delta_0 =
\emptyset)$.  We prove this by induction on $k$.

Let $k=0$.   The polynomial $G_{{\bf n},0} = A \Psi_{{\bf n},0}
+B\Psi_{{\bf n}_l,0}$ is not identically equal to zero, and $|{\bf
n}| \leq \deg(G_{{\bf n},0}) \leq |{\bf n}| +1$. Therefore, $
G_{{\bf n},0}$ has at most $|{\bf n}| +1$ zeros in $\mathbb{C}$.
Let $h_j, j=1,\ldots,m ,$ denote polynomials, where $\deg(h_j)
\leq n_j -1.$ According to (\ref{eq2}),
\begin{equation} \label{eq:a}
\int_{\Delta_1} G_{{\bf n},0}(x) \sum_{j=1}^m
h_j(x)\widehat{s}_{{2,j}}(x) d\sigma_1(x) = 0
\end{equation}
$(\widehat{s}_{{2,1}}\equiv 1)$.

In the sequel, we call change knot a point on the real line where
a function changes its sign. Notice that for each $k \in \{
0,\ldots,m-1\}$,   $G_{{\bf n},k}$ is a real function when
restricted to the real line. Assume that $G_{{\bf n},0}$ has $N
\leq |{\bf n}|-1$ change knots in the interior of $\Delta_1$. We
can find polynomials $h_j\,, j=1,\ldots,m, \deg(h_j) \leq n_j-1,$
such that $\sum_{j=1}^m h_j\widehat{s}_{{2,j}}$ has a simple zero
at each change knot of $G_{{\bf n},0}$ on $\Delta_1$ and a zero of
order $|{\bf n}|-1-N$ at one of the extreme points of $\Delta_1$.
By Lemma \ref{lm:AT},
$(1,\widehat{s}_{{2,2}},\ldots,\widehat{s}_{{2,m}})$ forms an AT
system with respect to ${\bf n} \in \mathbb{Z}_+^m(*)$; therefore,
$\sum_{j=1}^m h_j\widehat{s}_{{2,j}}$ can have no other zero on
$\Delta_1$, but this contradicts (\ref{eq:a}) since $G_{{\bf n},0}
\sum_{j=1}^m h_j \widehat{s}_{{2,j}} $ would have a constant sign
on $\Delta_1$ (and $\supp{\sigma_1}$ contains infinitely many
points). Therefore, $G_{{\bf n},0}$ has at least $|{\bf n}|$
change knots in the interior of  $\Delta_1$. Consequently, all the
zeros of $G_{{\bf n},0}$ are simple and lie on $\mathbb{R}$ as
claimed.

Assume that for each $k \in \{0,\ldots,\kappa-1\}, 1\leq \kappa
\leq m-1,$ the claim is satisfied whereas it is violated when
$k=\kappa$. Let $h_j$ denote polynomials such that $\deg(h_j) \leq
n_{j}^{\kappa} -1, \kappa +1 \leq j \leq m.$ Using (\ref{eq2}) or
(\ref{eq10})-(\ref{eq13}) according to the situation (to simplify
the writing we use the notation of (\ref{eq2}) but the arguments
are the same when $m=3$ and $n_1 < n_2 < n_3$; in particular, in
this case, $ds_{r_{0}}^{0} = ds_{1,3}, ds_{r_{1}}^{1} =
\widehat{s}_{3,2}d\tau_{2,3}/\widehat{\sigma}_3$ and
$ds_{r_{2}}^{2} = \widehat{s}_{2,3}d\tau_{3,2}/\widehat{\sigma}_2
$)
\begin{equation}
\label{eq:b} \int_{\Delta_{\kappa+1}} G_{{\bf n},\kappa}(x)
\sum_{j=\kappa+1}^m h_j(x) \widehat{s}_{\kappa+2,j}^{\kappa}(x)
d\sigma^{\kappa}_{\kappa+1}(x) = 0\,
\end{equation}
$(\widehat{s}_{\kappa+2,\kappa+1}^{\kappa}\equiv 1)$. Arguing as
above, since $(1,\widehat{s}_{\kappa+2,\kappa+2}^{\kappa},
\ldots,\widehat{s}_{\kappa+2,m}^{\kappa})$ forms an AT system with
respect to ${\bf n}^{\kappa} \in \mathbb{Z}_+^{m-\kappa}(*)$, we
conclude that $G_{{\bf n},\kappa}$ has at least $|{\bf
n}^{\kappa}|$ change knots in the interior of $\Delta_{\kappa+1}$.

Let us suppose that $G_{{\bf n},\kappa}$ has at least $|{\bf
n}^{\kappa}|+2$ zeros in ${\mathbb{C}}\setminus \Delta_{\kappa}$
and let $W_{{\bf n},\kappa}$ be the monic polynomial whose zeros
are those points (counting multiplicities). The complex zeros of
$G_{{\bf n},\kappa}$ (if any) must appear in conjugate pairs since
$G_{{\bf n},\kappa}(\overline{z}) = \overline{G_{{\bf n},\kappa}(
{z})};$ therefore, the coefficients of $W_{{\bf n},\kappa}$ are
real numbers. On the other hand, from (\ref{eq2}) ((\ref{eq10}) or
(\ref{eq12}) when necessary)
\[  0 =
\int_{\Delta_{\kappa}} G_{{\bf n},\kappa-1}(x)
\frac{z^{n^{\kappa-1}_{r_{\kappa-1}}} -
x^{n^{\kappa-1}_{r_{\kappa-1}}}}{z-x}
ds_{r_{\kappa-1}}^{\kappa-1}(x)\,.
\]
Therefore,
\[ G_{{\bf n},\kappa}(z) = \frac{1}{z^{n^{\kappa-1}_{r_{\kappa-1}}}} \int_{\Delta_{\kappa}}
\frac{x^{n^{\kappa-1}_{r_{\kappa-1}}}G_{{\bf n},\kappa-1}(x)}{z-x}
ds_{r_{\kappa-1}}^{\kappa-1}(x) = {\mathcal{O}}
\left(\frac{1}{z^{n^{\kappa-1}_{r_{\kappa-1}}+1}}\right)\,, \quad
z \to \infty\,,
\]
and taking into consideration the degree of $W_{{\bf
n},\kappa}\,,$ we obtain
\[ \frac{z^jG_{{\bf n},\kappa}}{W_{{\bf n},\kappa}} = {\mathcal{O}}
\left(\frac{1}{z^{2}}\right)\in H({\mathbb{C}} \setminus
\Delta_{\kappa})\,,\qquad j = 0,\ldots, |{\bf n}^{\kappa-1}| + 1
\,.
\]

Let $\Gamma$ be a closed Jordan curve which surrounds
$\Delta_{\kappa}$ and such that all the zeros of $W_{{\bf
n},\kappa}$ lie in the exterior of $\Gamma$. Using Cauchy's
Theorem, the integral expression for $G_{{\bf n},\kappa}$,
Fubini's Theorem, and Cauchy's Integral Formula, for each $j =
0,\ldots, |{\bf n}^{\kappa-1}|+ 1$, we have
\[ 0 = \frac{1}{2\pi i}\int_{\Gamma} \frac{z^j G_{{\bf n},\kappa}(z)}{W_{{\bf n},\kappa}(z)}
dz =  \frac{1}{2\pi i}\int_{\Gamma} \frac{z^j}{W_{{\bf
n},\kappa}(z)}\int_{\Delta_{\kappa}} \frac{G_{{\bf
n},\kappa-1}(x)}{z-x} ds_{r_{\kappa-1}}^{\kappa-1}(x) dz =
\]
\[ \int_{\Delta_{\kappa}} \frac{x^j G_{{\bf n},\kappa-1}(x)}{W_{{\bf
n},\kappa}(x)} ds_{r_{\kappa-1}}^{\kappa-1}(x)\,,
\]
which implies that $G_{{\bf n},\kappa-1}$ has at least $|{\bf
n}^{\kappa-1}|+ 2$ change knots in the interior of
$\Delta_{\kappa}$. This contradicts our induction hypothesis since
this function can have at most $|{\bf n}^{\kappa-1}|+ 1$ zeros in
${\mathbb{C}} \setminus \Delta_{\kappa-1} \supset
\Delta_{\kappa}$. Hence  $G_{{\bf n},\kappa}$ has at most $|{\bf
n}^{\kappa}|+1$ zeros in ${\mathbb{C}}\setminus \Delta_{\kappa}$
as claimed.

Taking $B=0$ the assumption ${\bf n}_l \in \mathbb{Z}_+^m(*,\tau)$
is not required, and the arguments above lead to the proof that
$\Psi_{{\bf n},k}$ has at most $|{\bf n}^k|$ zeros on $\mathbb{C}
\setminus \Delta_k$ since $Q_{{\bf n}}= \Psi_{{\bf n},0}$ has at
most $|{\bf n}|$ zeros on $\mathbb{C}$. Consequently, the zeros of
$\Psi_{{\bf n},k}$ in $\mathbb{C} \setminus \Delta_k$ are exactly
the  $|{\bf n}^k|$ simple ones it has in the interior of
$\Delta_{k+1}$.

Let $I$ be the closure of a connected component of $\Delta_{k+1}
\setminus \supp{\sigma_{k+1}^{k}}$ and let us assume that $I$
contains two consecutive simple zeros $x_1,x_2$ of $\Psi_{{\bf
n},k}$. Taking $B=0$ and $A=1$, we can rewrite (\ref{eq:b}) as
follows
\begin{equation}\label{eq:c}
\int_{\Delta_{k+1}} \frac{\Psi_{{\bf n},k}(x)}{(x-x_1)(x-x_2)}
\sum_{j=k+1}^m h_j(x) \widehat{s}_{k+2,j}^{k}(x)
(x-x_1)(x-x_2)d\sigma^{k}_{k+1}(x) = 0\,,
\end{equation}
where $\deg(h_j)\leq n_j^k-1, j=k+1,\ldots,m.$ The measure
$(x-x_1)(x-x_2)d\sigma^{k}_{k+1}(x)$ has a constant sign on
$\Delta_{k+1}$ and $ {\Psi_{{\bf n},k}(x)}/{(x-x_1)(x-x_2)}$ has
$|{\bf n}^k| -2$ change knots on $\Delta_{k+1}$. Using again Lemma
\ref{lm:AT}, we can construct appropriate polynomials $h_j$ to
contradict (\ref{eq:c}). Consequently, such an interval $I$ cannot
exist.

Fix $y \in {\mathbb{R}}\setminus \Delta_{k}$ and $k \in
\{0,1,\ldots,m-1\}$. It  cannot occur  that $\Psi_{{\bf n}_l,k}(y)
= \Psi_{{\bf n},k}(y) = 0$. If this was so, $y$ would have to be a
simple zero of $\Psi_{{\bf n}_l,k}$ and $\Psi_{{\bf n},k}$.
Therefore, $(\Psi_{{\bf n}_l,k})^{\prime}(y) \neq 0 \neq
(\Psi_{{\bf n},k})^{\prime}(y)$. Taking $A=1,B=-\Psi_{{\bf
n},k}^{\prime}(y)/\Psi_{{\bf n}_l,k}^{\prime}(y)$, we find that
\[ G_{{\bf n},k}(y)  = (A\Psi_{{\bf n},k}+B\Psi_{{\bf
n}_{l},k})(y) = (G_{{\bf n},k})^{\prime}(y) = 0\,,
\]
which means that $G_{{\bf n},k}$ has at least a double zero at $y$
against what we proved before.

Now, taking $A = \Psi_{{\bf n}_l,k}(y), B = -\Psi_{{\bf n},k}
(y)$, we have that $|A| + |B| > 0$. Since
\[ \Psi_{{\bf n}_l,k}(y)\Psi_{{\bf n},k}(y) - \Psi_{{\bf
n},k}(y)\Psi_{{\bf n}_l,k}(y) = 0 \,,
\]
and the zeros on $\mathbb{R} \setminus \Delta_{k}$ of $\Psi_{{\bf
n}_l,k}(y)\Psi_{{\bf n},k}(x) - \Psi_{{\bf n},k}(y)\Psi_{{\bf
n}_l,k}(x)$ with respect to $x$ are simple, using again what we
proved above, it follows that
\[ \Psi_{{\bf n}_l,k}(y)\Psi_{{\bf n},k}^{\prime}(y) - \Psi_{{\bf
n},k}(y)\Psi_{{\bf n}_l,k}^{\prime}(y) \neq 0\,.
\]
But $\Psi_{{\bf n}_l,k}(y)\Psi_{{\bf n},k}^{\prime}(y) -
\Psi_{{\bf n},k}(y)\Psi_{{\bf n}_l,k}^{\prime}(y)$ is a continuous
real function on $\mathbb{R} \setminus \Delta_k$ so it must have
constant sign on each one of the intervals forming $ \mathbb{R}
\setminus \Delta_{k}$; in particular, its sign on $\Delta_{k+1}$
is constant.

We know that $\Psi_{{\bf n}_l,k}$ has at least $|{\bf n}^k|$
simple zeros in the interior of $\Delta_{k+1}$. Evaluating
$\Psi_{{\bf n}_l,k}(y)\Psi_{{\bf n},k}^{\prime}(y) - \Psi_{{\bf
n},k}(y)\Psi_{{\bf n}_l,k}^{\prime}(y)$ at two consecutive zeros
of $\Psi_{{\bf n}_l,k}$, since the sign of $\Psi^{\prime}_{{\bf
n}_l,k}$ at these two points changes  the sign of $\Psi_{{\bf
n},k}$ must also change. Using Bolzano's theorem we find that
there must be an intermediate zero of $\Psi_{{\bf n},k}$.
Analogously, one proves that between two consecutive zeros of
$\Psi_{{\bf n} ,k}$ on $\Delta_{k+1}$ there is one of $\Psi_{{\bf
n}_l,k}$. Thus, the interlacing property has been proved.
\end{proof}

Let $Q_{{\bf n},k+1}, k=0,\ldots,m-1,$ denote the monic polynomial
whose zeros are equal to those of $\psink$ on $\Delta_{k+1}$. From
(\ref{eq2}) ((\ref{eq10}), (\ref{eq12}),  or (\ref{eq13})   when
necessary)
\[  0 =
\int_{\Delta_{k+1}} \Psi_{{\bf n},k}(x) \frac{z^{n^{k}_{r_{k}}} -
x^{n^{k }_{r_{k}}}}{z-x} ds_{r_{k}}^{k}(x)
\]
(Recall that when $m=3$ and $n_1 < n_2 < n_3$, we take
$ds_{r_{0}}^{0} = ds_{1,3}, ds_{r_{1}}^{1} =
\widehat{s}_{3,2}d\tau_{2,3}/\widehat{\sigma}_3$ and
$ds_{r_{2}}^{2} =
\widehat{s}_{2,3}d\tau_{3,2}/\widehat{\sigma}_2$.) Therefore,
\[ \Psi_{{\bf n},k+1}(z) = \frac{1}{z^{n^{k}_{r_{k}}}} \int_{\Delta_{k+1}}
\frac{x^{n^{k}_{r_{k}}}\Psi_{{\bf n},k}(x)}{z-x} ds_{r_{k}}^{k}(x)
= {\mathcal{O}} \left(\frac{1}{z^{n^{k}_{r_{k}}+1}}\right)\,,
\quad z \to \infty\,,
\]
and taking into consideration the degree of $Q_{{\bf n},k+2}$ (by
definition $Q_{{\bf n},m+1} \equiv 1$), we obtain
\[ \frac{z^j\Psi_{{\bf n},k+1}}{Q_{{\bf n},{k+2}}} = {\mathcal{O}}
\left(\frac{1}{z^{2}}\right)\in H({\mathbb{C}} \setminus
\Delta_{k+1})\,,\qquad j = 0,\ldots, |{\bf n}^{k}| - 1 \,.
\]

Let $\Gamma$ be a closed Jordan curve which surrounds
$\Delta_{k+1}$ such that all the zeros of $Q_{{\bf n},k+2}$ lie in
the exterior of $\Gamma$. Using Cauchy's Theorem, the integral
expression for $\Psi_{{\bf n},k+1}$, Fubini's Theorem, and
Cauchy's Integral Formula, for each $j = 0,\ldots, |{\bf n}^{k}|-
1$ (we also define $Q_{{\bf n},0} \equiv 1$), we have
\[ 0 = \frac{1}{2\pi i}\int_{\Gamma} \frac{z^j \Psi_{{\bf n},k+1}(z)}{Q_{{\bf n},k+2}(z)}
dz =  \frac{1}{2\pi i}\int_{\Gamma} \frac{z^j}{Q_{{\bf
n},k+2}(z)}\int_{\Delta_{k+1}} \frac{\Psi_{{\bf n},k}(x)}{z-x}
ds_{r_{k}}^{k}(x) dz =
\]
\begin{equation} \label{eq:orto1} \int_{\Delta_{k+1}} x^j Q_{{\bf n},k+1}(x)  \frac{H_{{\bf
n}, k+1}(x)ds_{r_{k}}^{k}(x)}{Q_{{\bf n},k}(x)Q_{{\bf
n},k+2}(x)}\,, \quad k=0,\ldots,m-1\,,
\end{equation}
where
\[ H_{{\bf n}, k+1} = \frac{Q_{{\bf n},k} \Psi_{{\bf n},k} }{Q_{{\bf
n},k+1} }\,, \qquad k=0,\ldots,m\,,
\]
has constant sign on $\Delta_{k+1}$.

This last relation implies  that
\[
 \int_{\Delta_{k+1}}  \frac{(Q(z) - Q(x))}{z-x}
\,\,Q_{{\bf n},k+1}(x)\,\frac{H_{{\bf n},
k+1}(x)ds_{r_{k}}^{k}(x)}{Q_{{\bf n},k}(x)Q_{{\bf n},k+2}(x)} =
0\,,
\]
where $Q$ is any polynomial of degree $\leq |{\bf n}^k|$. If we
use this formula with $Q = Q_{{\bf n},k+1}$ and $Q= Q_{{\bf
n},k+2}$, respectively, we obtain
\[
\int_{\Delta_{k+1}} \frac{Q_{{\bf n},k+1}(x)}{z-x}\frac{H_{{\bf
n}, k+1}(x)ds_{r_{k}}^{k}(x)}{Q_{{\bf n},k}(x)Q_{{\bf n},k+2}(x)}
=
\]
\[ \frac{1}{Q_{{\bf n},k+1}(z)}\int_{\Delta_{k+1}} \frac{ Q^2_{{\bf
n},k+1}(x)}{z-x}\frac{H_{{\bf n},
k+1}(x)ds_{r_{k}}^{k}(x)}{Q_{{\bf n},k}(x)Q_{{\bf n},k+2}(x)}
\]
and
\[
\int_{\Delta_{k+1}} \frac{Q_{{\bf n},k+1}(x)}{z-x}\frac{H_{{\bf
n}, k+1}(x)ds_{r_{k}}^{k}(x)}{Q_{{\bf n},k}(x)Q_{{\bf n},k+2}(x)}
=
\]
\[\frac{1}{Q_{{\bf n},k+2}(z)}\int_{\Delta_{k+1}} \frac{\Psi_{{\bf
n},k}(x)ds_{r_{k}}^{k}(x) }{z-x}  \,.
\]
Equating these two relations and using the definition of
$\Psi_{{\bf n},k+1}$ and $H_{{\bf n},k+2}$, we obtain
\begin{equation} \label{eq:5.7}
H_{{\bf n},k+2}(z)=\int_{\Delta_{k+1}} \frac{ Q^2_{{\bf
n},k+1}(x)}{z-x}\frac{H_{{\bf n},
k+1}(x)ds_{r_{k}}^{k}(x)}{Q_{{\bf n},k}(x)Q_{{\bf n},k+2}(x)}\,,
\quad k=0,\ldots,m-1\,.
\end{equation}
Notice that from the definition $H_{{\bf n},1}\equiv 1.$

For each $k=1,\ldots,m,$ set
\begin{equation} \label{eq:K} K_{{\bf n},k}^{-2} =   \int_{\Delta_k}
Q_{{\bf n},k}^2(x) \left| \frac{Q_{{\bf n},k-1}(x)\Psi_{{\bf
n},k-1}(x)}{Q_{{\bf n},k}(x)}\right|
\frac{d|s^{k-1}_{r_{k-1}}|(x)}{|Q_{{\bf n},k-1}(x)Q_{{\bf
n},k+1}(x)|}  \;,
\end{equation}
where $|s|$ denotes the total variation of the measures $s$. Take
\[ K_{{\bf n},0} = 1 \;, \quad \kappa_{{\bf n},k} = \frac{K_{{\bf n},k}}{K_{{\bf n},k-1}}
\;, \quad k=1,\ldots,m \;.
\]
Define
\begin{equation} \label{eq:orton} q_{{\bf n},k} =
\kappa_{{\bf n},k}Q_{{\bf n},k} \;, \quad h_{{\bf n},k}  = K_{{\bf
n},k-1}^2 H_{{\bf n},k} \;,
\end{equation}
and
\begin{equation} \label{eq:var1}
d\rho_{{\bf n},k}(x) = \frac{h_{{\bf
n},k}(x)ds^{k-1}_{r_{k-1}}(x)}{Q_{{\bf n},k-1}(x)Q_{{\bf
n},k+1}(x)} \,.
\end{equation}
Notice that the measure $\rho_{{\bf n},k}$ has constant sign on
$\Delta_k$. Let $\varepsilon_{{\bf n},k}$ be the sign of
$\rho_{{\bf n},k}$. From (\ref{eq:orto1}) and the notation
introduced above, we obtain
\begin{equation} \label{orton1}
\int_{\Delta_k} x^{\nu} q_{{\bf n},k}(x)d|\rho_{{\bf
n},k}|(x)=0,\quad \nu=0,\ldots,|{\bf n}^{k-1}| -1,\quad
k=1,\ldots,m\,,
\end{equation}
and $q_{n,k}$ is orthonormal with respect to the varying measure
$|\rho_{{\bf n},k}|.$ On the other hand, using (\ref{eq:5.7}) it
follows that
\begin{equation} \label{eq:5.7*}
h_{{\bf n},k+1}(z)= \varepsilon_{{\bf n},k}\int_{\Delta_{k}}
\frac{ q^2_{{\bf n},k}(x)}{z-x}d|\rho_{{\bf n},k}|(x)\,, \quad
k=1,\ldots,m\,.
\end{equation}

\begin{lem} \label{lm:3}
Let $S = {\mathcal{N}}(\sigma_1,\ldots,\sigma_m)$ be a Nikishin
system such that $ \supp{\sigma_k}  = \widetilde{\Delta}_k \cup
e_k, k=1,\ldots,m$, where $\widetilde{\Delta}_k$ is a bounded
interval of the real line, $\sigma_k^{\prime} > 0$ a.e. on
$\widetilde{\Delta}_k$, and $e_k$ is a set without accumulation
points in $\mathbb{R} \setminus \widetilde{\Delta}_k$. Let
$\Lambda \subset {\mathbb{Z}}_+^m(*)$ be an infinite sequence of
distinct multi-indices such that  $\displaystyle{\max_{{\bf n} \in
\Lambda }(\max_{k =1,\ldots,m}mn_{k} -  |{\bf n}| ) }< \infty.$
For any continuous function $f$ on $\supp{\sigma_{k}^{k-1}}$
\begin{equation} \label{eq:5.7**}
\lim_{{\bf n} \in \Lambda}  \int_{\Delta_{k}} f(x)q^2_{{\bf
n},k}(x) d|\rho_{{\bf n},k}|(x) = \frac{1}{\pi}
\int_{\widetilde{\Delta}_k} f(x) \frac{d x}{\sqrt{(b_k -x)(x -
a_k)}} \,,
\end{equation}
where $\widetilde{\Delta}_k =[a_k,b_k]$. In particular,
\begin{equation} \label{eq:h}
\lim_{{\bf n} \in \Lambda} \varepsilon_{{\bf n},k}h_{{\bf
n},k+1}(z) = \frac{1}{\sqrt{(z - b_k)(z - a_k)}}\,,\quad
\mathcal{K} \subset \mathbb{C} \setminus
\supp{\sigma_{k}^{k-1}}\,,
\end{equation}
where $\sqrt{(z - b_k)(z - a_k)} > 0$ if $z > 0$. Consequently,
for $k = 1,\ldots,m,$ each point of $\supp{\sigma^{k-1}_k}
\setminus \widetilde{\Delta}_k,$ is a limit of zeros of $\{Q_{{\bf
n},k}\}, {\bf n} \in \Lambda.$
\end{lem}
\begin{proof} We will proof this by induction on $k$. For $k=1$,
using Corollary 3 in \cite{BCL}, it follows that
\[ \lim_{{\bf n} \in \Lambda}  \int_{\Delta_{1}} f(x)q^2_{{\bf
n},1}(x) \frac{ d|s_{r_{0}}^{0}|(x)}{| Q_{{\bf n},2}(x)|} =
\frac{1}{\pi} \int_{\widetilde{\Delta}_1} f(x) \frac{d
x}{\sqrt{(b_1 -x)(x - a_1)}}\,,
\]
where $f$ is continuous on $\supp{\sigma_1}$. Take $f(x) =
(z-x)^{-1}$ where $z \in \mathbb{C} \setminus \supp{\sigma_1}$.
According to (\ref{eq:5.7*}) and the previous limit one obtains
that
\[ \lim_{{\bf n} \in \Lambda} \varepsilon_{{\bf n},1}h_{{\bf n},2}(z) =
\frac{1}{\sqrt{(z - b_1)(z - a_1)}} =: h_2(z)\,,
\]
pointwise on $\mathbb{C} \setminus \supp{\sigma_1}$. Since
\[\left| \int_{\Delta_{1}} \frac{q^2_{{\bf
n},1}(x)}{z-x} \frac{ d|s_{r_{0}}^{0}|(x)}{|Q_{{\bf
n},2}(x)|}\right| \leq \frac{1}{d(\mathcal{K},\supp{\sigma_1})}\,,
\quad z \in \mathcal{K} \subset \mathbb{C} \setminus
\supp{\sigma_1}\,,
\]
where $d(\mathcal{K},\supp{\sigma_1})$ denotes the distance
between the two compact sets, the sequence $\{h_{{\bf n},2}\},
{\bf n} \in \Lambda,$ is uniformly bounded on compact subsets of
$\mathbb{C} \setminus \supp{\sigma_1}$ and (\ref{eq:h}) follows
for $k=1$.

Let $\zeta \in \supp{\sigma_1} \setminus \widetilde{\Delta}_1$.
Take $r > 0$ sufficiently small so that the circle $C_{r} = \{z:
|z - \zeta| = r \}$ surrounds no other point of $\supp{\sigma_1}
\setminus \widetilde{\Delta}_1$ and contains no zero of $q_{{\bf
n},1}, {\bf n} \in \Lambda$. From (\ref{eq:h}) for $k=1$
\[ \lim_{{\bf n} \in \Lambda} \frac{1}{2\pi i}\int_{C_r}
\frac{\varepsilon_{{\bf n},1}h_{{\bf
n},2}^{\prime}(z)}{\varepsilon_{{\bf n},1}h_{{\bf n},2}(z)} d z=
\frac{1}{2\pi i}\int_{C_r} \frac{h_2^{\prime}(z)}{h_2(z)}dz = 0\,.
\]
From the definition, $\Psi_{{\bf n },1}, {\bf n} \in \Lambda,$ has
either a simple pole at $\zeta$ or $Q_{{\bf n},1}$ has a zero at
$\zeta$. In the second case there is nothing to prove. Let us
restrict our attention to those ${\bf n} \in \Lambda$ such that
$\Psi_{{\bf n },1}, {\bf n} \in \Lambda,$ has  a simple pole at
$\zeta$. Then, $h_{{\bf n},2} = K_{{\bf n},1}^2 Q_{{\bf
n},1}\Psi_{{\bf n},1}/Q_{{\bf n},2}$ also has a simple pole at
$\zeta$. Using the argument principle, it follows that for all
sufficiently large $|{\bf n}|, {\bf n} \in \Lambda$, $Q_{{\bf
n},1}$ must have a simple zero inside $C_r$. The parameter $r$ can
be taken arbitrarily small; therefore, the last statement of the
lemma readily follows and the basis of induction is fulfilled.

Let us assume that the lemma is satisfied for $k \in
\{1,\ldots,\kappa -1\}, 1\leq \kappa \leq m,$ and let us prove
that it is also true for $\kappa.$ From (\ref{eq:h}) applied to
$\kappa -1$, we have that
\[ \lim_{{\bf n} \in \Lambda} |h_{{\bf n},\kappa}(x)| =
\frac{1}{\sqrt{|(x - b_{\kappa -1)}(x - a_{\kappa -1})|}}\,,
\]
uniformly on $\Delta_{\kappa} \subset \mathbb{C} \setminus
\supp{\sigma_{\kappa -1}^{\kappa -2}}.$ It follows that
$\{|h_{{\bf n},\kappa}|d|s_{r_{\kappa-1}}^{\kappa-1}|\},{\bf n}
\in \Lambda,$ is a sequence of Denisov type measures according to
Definition 3 in \cite{BCL} and $(\{|h_{{\bf
n},\kappa}|d|s_{r_{\kappa-1}}^{\kappa-1}|\},\{|Q_{{\bf
n},\kappa-1} Q_{{\bf n},\kappa+1}|\},l),{\bf n} \in \Lambda,$ is
strongly admissible as in Definition 2 of \cite{BCL} for each $l
\in \mathbb{Z}$ (see paragraph just after both definitions in the
referred paper). Therefore, we can apply Corollary 3 in \cite{BCL}
of which (\ref{eq:5.7**}) is a particular case. In the proof of
Corollary 3 of \cite{BCL} (see also Theorem 9 in
\cite{BernardoGuillermo1}) it is required that $\deg (Q_{{\bf
n},k-1}Q_{{\bf n},k+1}) - 2\deg(Q_{{\bf n},k}) \leq C$ where $C
\geq 0$ is a constant. For $k=1$ this is trivially true (with
$C=0$). Since we apply an induction procedure on $k$, in order
that this requirement be fulfilled for all $k \in \{1,\ldots,m\}$
we impose that $\displaystyle{\max_{{\bf n} \in \Lambda }(\max_{k
=1,\ldots,m}mn_{k} -  |{\bf n}| ) }< \infty$. From
(\ref{eq:5.7**}), (\ref{eq:h}) and the rest of the statements of
the lemma immediately follow just as in the case when $k=1$. With
this we conclude the proof.
\end{proof}

\begin{rmk}
The last statement of Lemma \ref{lm:3} concerning the convergence
of the zeros of $Q_{{\bf n},1}$ outside $\widetilde{\Delta}_1$ to
the mass points of $\sigma_1$ on $\supp{\sigma_1} \setminus
\widetilde{\Delta}_1$ can be proved without the assumption that
$\sigma_k^{\prime} > 0$ a.e. on $\widetilde{\Delta}_k, k =
1,\ldots,m$. This is an easy consequence of Theorem 1 in
\cite{FidLop3}. From the proof of Lemma \ref{lm:3} it also follows
that if we only have $\sigma_k^{\prime} > 0$ a.e. on
$\widetilde{\Delta}_k, k = 1,\ldots,m^{\prime}, m^{\prime} \leq
m$, then (\ref{eq:5.7**})-(\ref{eq:h}) are satisfied for $k =
1,\ldots,m^{\prime}$ and the statement concerning the zeros holds
for $k=1,\ldots,m^{\prime}+1.$
\end{rmk}

\section{{\large Proof of main results}}

In this final section,  $S =
{\mathcal{N}}(\sigma_1,\ldots,\sigma_m)$ is a Nikishin system with
$ \supp{\sigma_k}  = \widetilde{\Delta}_k \cup e_k, k=1,\ldots,m$,
where $\widetilde{\Delta}_k$ is a bounded interval of the real
line, $\sigma_k^{\prime} > 0$ a.e. on $\widetilde{\Delta}_k$, and
$e_k$ is a set without accumulation points  in $\mathbb{R}
\setminus \widetilde{\Delta}_k$. Let $\Lambda \subset
\mathbb{Z}_{+}^{m}(*)$ be a sequence of distinct multi-indices.
Let us assume that there exists $l\in \{1,\ldots,m\}$ and a fixed
permutation $\tau$ of $\{1,\ldots,m\}$ such that for all ${\bf n}
\in \Lambda$ we have that ${\bf n},{\bf n}_l \in
\mathbb{Z}_{+}^{m}(*,\tau)$. From the interlacing property of the
zeros of $Q_{{\bf n},k}$ and $Q_{{\bf n}_l,k}$, and the limit
behavior of the zeros of $Q_{{\bf n},k}$ outside
$\widetilde{\Delta}_k$, it follows that the sequences
\[\left\{ {Q_{{\bf n}_l,k}}/{Q_{{\bf n},k}}\right\}_{{\bf n}\in\Lambda},\qquad
k=1,\ldots,m,
\]
are uniformly bounded on each compact subset of
$\mathbb{C}\setminus\supp{\sigma_{k}^{k-1}}$ for all sufficiently
large $|{\bf n}|$. By Montel's theorem, there exists a subsequence
of multi-indices $\Lambda' \subset \Lambda$ and a collection of
functions $\widetilde{F}^{l}_{k}$, holomorphic in
$\mathbb{C}\setminus\supp{\sigma_{k}^{k-1}}$, respectively, such
that
\begin{equation}\label{eq26}
\lim_{{\bf n}\in\Lambda'}\,\frac{Q_{{\bf n}_l,k}(z)}{Q_{{\bf
n},k}(z)}=\widetilde{F}^{(l)}_{k}(z),\quad \mathcal{K}
\subset\mathbb{C}\setminus\supp{\sigma_k^{k-1}},\,\,k=1,\ldots,m.
\end{equation}

In principle, the functions $\widetilde{F}^{(l)}_{k}$ may depend
on $\Lambda'$. We shall see that this is not the case and,
therefore, the limit in (\ref{eq26}) holds for ${\bf n} \in
\Lambda$. First, let us obtain some general information on the
functions $\widetilde{F}^{(l)}_{k}$.

The points in $\supp{\sigma^{k-1}_k} \setminus
\widetilde{\Delta}_k$ are isolated singularities of
$\widetilde{F}^{(l)}_{k}$. Let $\zeta \in \supp{\sigma^{k-1}_k}
\setminus \widetilde{\Delta}_k$. By Lemma \ref{lm:3} each such
point is a limit of zeros of $Q_{{\bf n},k}$ and $Q_{{\bf n}_l,k}$
as $|{\bf n}| \to \infty, {\bf n} \in \Lambda,$ and in a
sufficiently small neighborhood of them, for each ${\bf n} \in
\Lambda$, there can be at most one such zero of these polynomials
(so there is exactly one, for all sufficiently large $|{\bf n}|$).
Let $\lim_{{\bf n} \in \Lambda} \zeta_{\bf n} = \zeta$ where
$Q_{{\bf n},k}(\zeta_{\bf n}) =0$. From (\ref{eq26})
\[ \lim_{{\bf n}\in\Lambda'}\,\frac{(z - \zeta_{\bf n})Q_{{\bf n}_l,k}(z)}{Q_{{\bf
n},k}(z)}=(z - \zeta)\widetilde{F}^{(l)}_{k}(z),\quad \mathcal{K}
\subset (\mathbb{C}\setminus\supp{\sigma_k^{k-1}}) \cup
\{\zeta\}\,,
\]
and  $(z - \zeta)\widetilde{F}^{(l)}_{k}(z)$ is analytic in a
neighborhood of $\zeta$. Hence $\zeta$ is not an essential
singularity of $\widetilde{F}^{(l)}_{k}$. Taking into
consideration that $Q_{{\bf n}_l,k}, {\bf n}\in \Lambda$ also has
a sequence of zeros converging to $\zeta$, from the argument
principle it follows that $\zeta$ is a removable singularity of
$\widetilde{F}^{(l)}_{k}$ which is not a zero. Using the
interlacing property and the convergence of the zeros of $Q_{{\bf
n},k}$ and $Q_{{\bf n}_l,k}$ outside $\widetilde{\Delta}_k$ as
$|{\bf n}| \to \infty, {\bf n} \in \Lambda,$ to the points in
$\supp{\sigma_{k}^{k-1}} \setminus \widetilde{\Delta}_k$, it is
easy to deduce that on each compact subset of
$\mathbb{C}\setminus\supp{\sigma_k^{k-1}}$ the functions $|Q_{{\bf
n}_l,k}/Q_{{\bf n},k}|, {\bf n} \in \Lambda,$ are uniformly
bounded from below by a positive constant for all sufficiently
large $|{\bf n}|$. Therefore, in
$\mathbb{C}\setminus\supp{\sigma_k^{k-1}}$ the function
$\widetilde{F}^{(l)}_{k}$ is also different from zero. According
to the definition of $Q_{{\bf n},k}$ and $Q_{{\bf n}_l,k}$ and
Lemma \ref{entrelazamiento}, for $k=1,\ldots,\tau^{-1}(l)$, we
have that $\deg{Q_{{\bf n}_l,k}} = |{\bf n}_l^{k-1}| = |{\bf
n}^{k-1}|+1= \deg{Q_{{\bf n},k}} +1$ whereas, for
$k=\tau^{-1}(l)+1,\ldots,m$, we obtain that $\deg{Q_{{\bf n}_l,k}}
= |{\bf n}_l^{k-1}| = |{\bf n}^{k-1}|= \deg{Q_{{\bf n},k}}$.
Consequently, for $k=1,\ldots,\tau^{-1}(l)$, the function
$\widetilde{F}^{(l)}_{k}$ has a simple pole at infinity and
$(\widetilde{F}^{(l)}_{k})^{\prime}(\infty) =1$, whereas, for
$k=\tau^{-1}(l)+1,\ldots,m$, it is analytic at infinity and
$\widetilde{F}^{(l)}_{k}(\infty) =1$.

Now let us prove that the functions $\widetilde{F}^{(l)}_{k}$
satisfy a system of boundary value problems.

\begin{lem} \label{sistema}
Let $S = {\mathcal{N}}(\sigma_1,\ldots,\sigma_m)$ be a Nikishin
system with $ \supp{\sigma_k}  = \widetilde{\Delta}_k \cup e_k,
k=1,\ldots,m$, where $\widetilde{\Delta}_k$ is a bounded interval
of the real line, $\sigma_k^{\prime} > 0$ a.e. on
$\widetilde{\Delta}_k$, and $e_k$ is a set without accumulation
points in $\mathbb{R} \setminus \widetilde{\Delta}_k$. Let
$\Lambda \subset \mathbb{Z}_{+}^{m}(*)$ be a sequence of distinct
multi-indices such that $\displaystyle{\max_{{\bf n} \in \Lambda
}(\max_{k =1,\ldots,m}mn_{k} -  |{\bf n}| ) }< \infty.$ Let us
assume that there exists $l\in \{1,\ldots,m\}$ and a fixed
permutation $\tau$ of $\{1,\ldots,m\}$ such that for all ${\bf n}
\in \Lambda$ we have that ${\bf n},{\bf n}_l \in
\mathbb{Z}_{+}^{m}(*,\tau)$. Take $\Lambda'\subset\Lambda$ such
that $(\ref{eq26})$ holds. Then, there exists a normalization
$F^{(l)}_{k}$, $k=1,\ldots,m,$ by positive constants, of the
functions $\widetilde{F}^{(l)}_{k}$, $k=1,\ldots,m,$ given in
$(\ref{eq26})$, which verifies the system of boundary value
problems
\begin{equation} \label{sisec}
\aligned &1) \qquad
 F^{(l)}_{k},\, {1}/{F^{(l)}_{k}}\in
H(\mathbb{C}\setminus \widetilde{\Delta}_{k})\,,\\
&2) \qquad (F^{(l)}_{k})'(\infty)>0\,,\quad k=1,\ldots,\tau^{-1}(l)\,,\\
&2') \qquad F^{(l)}_{k}(\infty)>0\,,
\quad k=\tau^{-1}(l)+1,\ldots,m\,,\\
&3) \qquad
|F^{(l)}_{k}(x)|^{2}\frac{1}{|(F^{(l)}_{k-1}\,F^{(l)}_{k+1})(x)|}=1,\,\,
x\in \widetilde{\Delta}_{k}\,,\\
\endaligned
\end{equation}
 where $F^{(l)}_{0}\equiv\,F^{(l)}_{m+1}\equiv\,1$.
\end{lem}
\begin{proof}
The assertions 1), 2), and 2') were proved above for the functions
$\widetilde{F}_k^{(l)}$. Consequently, they are satisfied for any
normalization of these functions by means of positive constants.

From (\ref{orton1}) applied to ${\bf n}$ and ${\bf n}_l$, for each
$k=1,\ldots,m$, we have
\[
\int_{\Delta_k} x^{\nu} Q_{{\bf n},k}(x)d|\rho_{{\bf
n},k}|(x)=0,\qquad \nu=0,\ldots,|{\bf n}^{k-1}| -1\,,
\]
and
\[
\int_{\Delta_k} x^{\nu} Q_{{\bf n}_l,k}(x) g_{{\bf
n},k}(x)d|\rho_{{\bf n},k}|(x)=0\,, \qquad \nu=0,\ldots,|{\bf
n}_l^{k-1}| -1\,,
\]
where
\[ g_{{\bf n},k}(x) = \frac{|Q_{{\bf
n},k-1}(x)Q_{{\bf n},k+1}(x)|}{|Q_{{\bf n}_l,k-1}(x)Q_{{\bf
n}_l,k+1}(x)|}\frac{|h_{{\bf n}_l,k}(x)|}{|h_{{\bf n},k}(x)|}\,,
\quad d\rho_{{\bf n},k}(x) = \frac{h_{{\bf
n},k}(x)ds^{k-1}_{r_{k-1}}(x)}{Q_{{\bf n},k-1}(x)Q_{{\bf
n},k+1}(x)}\,.
\]
From (\ref{eq:h}) and (\ref{eq26})
\begin{equation} \label{eq:var}
\lim_{{\bf n} \in \Lambda'} g_{{\bf n},k}(x) =
 {|(\widetilde{F}_{k-1}^{(l)}\widetilde{F}_{k+1}^{(l)})(x)|^{-1}}
\end{equation}
uniformly on $\Delta_k$.

Fix $k \in \{\tau^{-1}(l)+1,\ldots,m\}$. As mentioned above, for
this selection of $k$ we have that $\deg Q_{{\bf n}_l,k} = \deg
Q_{{\bf n},k} = |{\bf n}^{k-1}|$. Using (3) in Theorem 1 and
Theorem 2 of \cite{BCL}, and (\ref{eq26}), it follows that
\begin{equation} \label{eq:otraa}
\lim_{{\bf n}\in \Lambda^{\prime}} \frac{Q_{{\bf
n}_l,k}(z)}{Q_{{\bf n},k}(z)}=\frac{S_k(z)}{S_k(\infty)} =
\widetilde{S}_k(z) = \widetilde{F}_k^{(l)}(z)\,, \qquad
\mathcal{K} \subset \overline{\mathbb{C}}\setminus
\supp{\sigma_k^{k-1}} \,,
\end{equation}
where $S_k$ denotes  the Szeg\H{o} function on
$\overline{\mathbb{C}} \setminus \widetilde{\Delta}_k$ with
respect to the weight
${|\widetilde{F}_{k-1}^{(l)}(x)\widetilde{F}_{k+1}^{(l)}(x)|^{-1}},
x \in \widetilde{\Delta}_k.$  The function $S_k$ is uniquely
determined by
\begin{equation}       \label{eq:Se}
\aligned &1)\qquad S_k,1/S_k \in H(\overline{\mathbb{C}} \setminus
\widetilde{\Delta}_k)\,,
\\
&2)\qquad  S_k(\infty)>0\,,
\\
&3 )\qquad |S_k(x)|^2
\frac1{\bigl|(\widetilde{F}_{k-1}^{(l)}\widetilde{F}_{k+1}^{(l)})(x)\bigr|}=1,
\qquad x \in \widetilde{\Delta}_k \,.
\endaligned
\end{equation}

Now, fix $k \in \{1, \dots , \tau^{-1}(l)\}$. In this situation
$\deg Q_{{\bf n}_l,k} = \deg Q_{{\bf n},k} +1 = |{\bf n}^{k-1}|
+1.$ Let $Q_{{\bf n},k}^*(x)$ be the monic polynomial of degree
$|{\bf n}^{k-1}|$ orthogonal with respect to the varying measure
$g_{{\bf n},k}d|\rho_{{\bf n},k}|$. Using the same arguments as
above, we have
\begin{equation}          \label{eq:comparativa}
\lim_{{\bf n}\in \Lambda^{\prime}} \frac{Q_{{\bf
n},k}^*(z)}{Q_{{\bf n},k}(z)}=\frac{S_k(z)}{S_k(\infty)} =
\widetilde{S}_k(z)\,, \qquad \mathcal{K} \subset
\overline{\mathbb{C}}\setminus \supp{\sigma_{k}^{k-1}} \,.
\end{equation}
On the other hand, since $\deg Q_{{\bf n}_l,k}=\deg Q_{{\bf
n},k}^* +1$ and both of these polynomials are orthogonal with
respect to the same varying weight,  by (3) and (4) in Theorem 1
of \cite{BCL} and (\ref{eq26}), it follows that
\begin{equation}                \label{eq:ratio}
\lim_{{\bf n}\in {\Lambda^{\prime}}} \frac{Q_{{\bf
n}_l,k}(z)}{Q_{{\bf n},k}^*(z)}=
\frac{\varphi_k(z)}{{\varphi_k^{\prime}}(\infty)} =
\widetilde{\varphi}_k(z)\,, \qquad \mathcal{K} \subset
\mathbb{C}\setminus \supp{\sigma_{k}^{k-1}}\,,
\end{equation}
where $\varphi_k$ denotes the conformal representation of
$\overline{\mathbb{C}}\setminus \widetilde{\Delta}_k$ onto $\{w :
|w|
> 1\}$ such that $\varphi_k(\infty) = \infty$ and
$\varphi_k^{\prime}(\infty)
> 0$. The function $\varphi_k$ is uniquely determined by
\begin{equation}    \label{eq:Zh}
   \aligned &1 )\qquad \varphi_k,1/\varphi_k \in
H(\mathbb{C} \setminus \widetilde{\Delta}_k)\,,
\\
&2)\qquad  \varphi_k^{\prime}(\infty)>0\,,
\\
&3 )\qquad |\varphi_k(x)|=1, \quad  x \in \widetilde{\Delta}_k \,.
\endaligned
\end{equation}
From (\ref{eq:comparativa}) and (\ref{eq:ratio}), we obtain
\begin{equation} \label{eq:otra}
\lim_{{\bf n}\in \Lambda^{\prime}} \frac{Q_{{\bf
n}_l,k}(z)}{Q_{{\bf n},k}(z)}=
(\widetilde{S}_k\widetilde{\varphi}_k)(z) =
\widetilde{F}_k^{(l)}(z)\,, \quad \mathcal{K} \subset
\mathbb{C}\setminus \supp{\sigma_k^{k-1}}  \,.
\end{equation}

Thus,
\begin{equation}    \label{eq:FS}
\widetilde{F}_k^{(l)}= \left\{
\begin{array}{ll} \widetilde{S}_k \widetilde{\varphi}_k\,, &  k = 1, \dots , \tau^{-1}(l)\,, \\
\widetilde{S}_k \,, &  k = \tau^{-1}(l)+1, \dots , m\,,
\end{array}
\right.
\end{equation}
and  from (\ref{eq:Se}) and (\ref{eq:FS}) it follows that
\begin{equation} \label{eq:sysinter}
|\widetilde{F}_k^{(l)}(x)|^2
\frac1{\bigl|(\widetilde{F}_{k-1}^{(l)}\widetilde{F}_{k+1}^{(l)})(x)\bigr|}=
\frac{1}{\omega_k}
 \,, \qquad x \in \widetilde{\Delta}_k\,, \qquad k=1,\ldots,m\,,
\end{equation}
where
\begin{equation} \label{eq:xd} \omega_k = \left\{
\begin{array}{ll}
(S_k\varphi_k^{\prime})^2(\infty)\,, & k=1,\ldots,\tau^{-1}(l)\,, \\
S_k^2(\infty)\,, & k=\tau^{-1}(l)+1,\ldots,m \,.
\end{array}
\right.
\end{equation}

Now, let us show that there exist positive constants $c_k,
k=1,\ldots,m,$ such that the functions ${F_k^{(l)}} = c_k
\widetilde{F}_k^{(l)}$ satisfy (\ref{sisec}). In fact, according
to (\ref{eq:sysinter}) for any such constants $c_k$ we have that
\[ |{F_k^{(l)}}(x)|^2
\frac1{\bigl|({F_{k-1}^{(l)}}{F_{k+1}^{(l)}})(x)\bigr|}=
\frac{c_k^2}{c_{k-1}c_{k+1}\omega_k}
 \,, \qquad x \in \widetilde{\Delta}_k\,,\qquad k=1,\ldots,m\,,
\]
where $c_0 = c_{m+1} = 1$. The problem reduces to finding
appropriate constants $c_k$ such that
\begin{equation} \label{eq:aa}
\frac{c_k^2}{c_{k-1}c_{k+1}\omega_k} = 1\,, \qquad k
=1,\ldots,m\,.
\end{equation}
Taking logarithm, we obtain the linear system of equations
\begin{equation} \label{eq:bb} 2\log c_k - \log c_{k-1} -
\log c_{k+1} = \log \omega_k \,, \qquad k=1,\ldots,m
\end{equation}
$(c_0 = c_{m+1} = 1)$ on the unknowns $\log c_k\,.$ This system
has a unique solution with which we conclude the proof.
\end{proof}

Consider the $(m+1)$-sheeted compact Riemann surface $\mathcal{R}$
introduced in Section 1.  Given $l \in\{1,\ldots,m\}$, let
$\psi^{(l)}$ be a singled valued function defined on $\mathcal{R}$
onto the extended complex plane satisfying
\[ \psi^{(l)}(z)=\frac{C_{1}}{z}+\mathcal{O}(\frac{1}{z^{2}}),\quad z\rightarrow\infty ^{(0)}\]
\[ \psi^{(l)}(z)=C_{2}\,z+\mathcal{O}(1),\quad z\rightarrow\infty ^{(l)}\]
where $C_{1}$ and $C_{2}$ are nonzero constants. Since the genus
of $\mathcal{R}$ is zero, $\psi^{(l)}$ exists and is uniquely
determined up to a multiplicative constant. Consider the branches
of $\psi^{(l)}$, corresponding to the different sheets
$k=0,1,\ldots,m$ of $\mathcal{R}$
\[\psi^{(l)}:=\{\psi^{(l)}_{k}\}_{k=0}^{m}\,. \]
We normalize $\psi^{(l)}$ so that
\[\prod_{k=0}^{m}\,\psi^{(l)}_{k}(\infty)=1\,. \]
Since the product of all the branches $\prod_{k=0}^{m}
\psi^{(l)}_{k}$ is a single valued analytic function in
$\overline{\mathbb{C}}$ without singularities, by Liouville's
Theorem it is constant and because of the normalization introduced
above
\[ \prod_{k=0}^{m}
\psi^{(l)}_{k}(z) \equiv 1\,, \qquad z \in
\overline{\mathbb{C}}\,.
\]

Given an arbitrary function $F(z)$  which has in a neighborhood of
infinity a Laurent expansion of the form $F(z)= Cz^k +
{\mathcal{O}}(z^{k-1}), C \neq 0,$ and $ k \in {\mathbb{Z}},$ we
denote
\begin{equation}            \label{eq:tild}
\widetilde{F}:= {F}/{C}\,.
\end{equation}
In particular,
\begin{equation} \label{eq:G} G_0^{(l)}(z) =
1/\widetilde{\psi}_0^{(l)}(z) = \prod_{k=1}^{m}
\widetilde{\psi}^{(l)}_{k}(z)\,.
\end{equation}

We are ready to state and prove our main result.

\begin{theo} \label{teofund}
Let $S = {\mathcal{N}}(\sigma_1,\ldots,\sigma_m)$ be a Nikishin
system with $ \supp{\sigma_k}  = \widetilde{\Delta}_k \cup e_k,
k=1,\ldots,m$, where $\widetilde{\Delta}_k$ is a bounded interval
of the real line, $\sigma_k^{\prime} > 0$ a.e. on
$\widetilde{\Delta}_k$, and $e_k$ is a set without accumulation
points in $\mathbb{R} \setminus \widetilde{\Delta}_k$. Let
$\Lambda \subset \mathbb{Z}_{+}^{m}(*)$ be a sequence of distinct
multi-indices such that  $\displaystyle{\max_{{\bf n} \in \Lambda
}(\max_{k =1,\ldots,m}mn_{k} -  |{\bf n}| ) }< \infty$. Let us
assume that there exists $l\in \{1,\ldots,m\}$ and a fixed
permutation $\tau$ of $\{1,\ldots,m\}$ such that for all ${\bf n}
\in \Lambda$ we have that ${\bf n},{\bf n}_l \in
\mathbb{Z}_{+}^{m}(*,\tau)$. Let $\{Q_{{\bf n},k}\}_{k=1}^{m}$,
${\bf n}\in\Lambda$, be the corresponding sequences of polynomials
defined in section 3. Then, for each fixed $k\in\{1,\ldots,m\}$,
we have
\begin{equation}\label{eq30}
\lim_{{\bf n}\in\Lambda}\,\frac{Q_{{\bf n}_l,k}(z)}{Q_{{\bf
n},k}(z)}=\widetilde{F}^{(l)}_{k}(z),\qquad z\in \mathcal{K}
\subset\mathbb{C}\setminus\supp{\sigma_k^{k-1}}
\end{equation}
where
\begin{equation}\label{eq31}
F^{(l)}_{k}:=\prod_{\nu=k}^{m}\,\psi^{(l)}_{\nu}.
\end{equation}
\end{theo}
\begin{proof} Since the families of functions
\[ \left\{ {Q_{{\bf n}_l,k}}/{Q_{{\bf n},k}}\right\}_{{\bf n}\in \Lambda}
\,,\qquad k=1,\ldots,m,
\]
are uniformly bounded on each compact subset $\mathcal{K} \subset
\mathbb{C} \setminus \supp{\sigma_k^{k-1}}$ for all sufficiently
large $|{\bf n}|, {\bf n} \in \Lambda$,  uniform convergence on
compact subsets of the indicated region follows from proving that
any system of convergent subsequences has the same limits.
According to Lemma \ref{sistema} the limit functions of such
convergent subsequences appropriately normalized always satisfy
the same system of boundary value problems (\ref{sisec}).
According to Lemma 4.2 in \cite{AptLopRoc} this system has a
unique solution and it is given by (\ref{eq31}). Since the
polynomials $Q_{{\bf n},k}$ and $Q_{{\bf n}_l,k}$ are monic, the
limit in (\ref{eq30}) must be the result of applying the action
$\widetilde{\mbox{}}\,\,$ defined in (\ref{eq:tild}) to
(\ref{eq31}).
\end{proof}

Theorem \ref{teo2} is a particular case of Theorem \ref{teofund}
on account of (\ref{eq:G}).

\medskip

\noindent {\sl Proof of Corollary \ref{cor1}.} Let
\[\Lambda_{\tau} =  \Lambda \cap \mathbb{Z}_+^m(*,\tau) \,,
\]
where $\tau$ is a given permutation of $\{1,\ldots,m\}$. We are
only interested in those $\Lambda_{\tau}$ with infinitely many
terms. There are at most $m!$ such subsequences. For ${\bf n} \in
\Lambda_{\tau}$ fixed, denote ${\bf n}_{\tau(j)}, j\in \{
1,\ldots,m\},$ the multi-index obtained adding one to all $j$
components $\tau(1),\ldots,\tau(j)$ of ${\bf n}$. (Notice that
this notation differs from that introduced previously for ${\bf
n}_l$.) In particular, ${\bf n}_{\tau(m)} = {\bf n+1}$.  It is
easy to verify that for all $j \in \{1,\ldots,m\}$, ${\bf
n}_{\tau(j)} \in \Lambda_{\tau}$. For all ${\bf n} \in
\Lambda_{\tau}$ and each $k \in \{1,\ldots,m\}$, we have
\[ \frac{Q_{{\bf n+1},k}}{Q_{{\bf n},k}} = \prod_{j=0}^{m-1}
\frac{Q_{{\bf n}_{\tau(j+1)},k}}{Q_{{\bf n}_{\tau(j)},k}}\,,
\]
where $Q_{{\bf n}_{\tau(0)},k} = Q_{{\bf n},k}$. From (\ref{eq30})
it follows that
\[ \lim_{{\bf n} \in \Lambda_{\tau}} \frac{Q_{{\bf n+1},k}(z)}{Q_{{\bf
n},k}(z)}= \prod_{l=1}^{m} \widetilde{F}_k^{(l)}(z)\,, \qquad
\mathcal{K} \subset \mathbb{C} \setminus \supp{\sigma_k^{k-1}}\,.
\]
The right side does not depend on $l$, since all possible values
intervene. Therefore, the limit is the same for all $\tau$ and
thus
\begin{equation} \label{ultima}
\lim_{{\bf n} \in \Lambda} \frac{Q_{{\bf n+1},k}(z)}{Q_{{\bf
n},k}(z)}= \prod_{l=1}^{m} \widetilde{F}_k^{(l)}(z)\,, \qquad
\mathcal{K} \subset \mathbb{C} \setminus \supp{\sigma_k^{k-1}}\,.
\end{equation}
Formula (\ref{eq:xe1}) is (\ref{ultima}) for $k=1$ on account of
(\ref{eq:G}) and (\ref{eq31}). \hfill $\Box$

The following corollary complements Theorem \ref{teofund}. The
proof is similar to that of Corollary 4.1 in \cite{AptLopRoc}.

\begin{co} \label{cor3}
Let $S = {\mathcal{N}}(\sigma_1,\ldots,\sigma_m)$ be a Nikishin
system with $ \supp{\sigma_k}  = \widetilde{\Delta}_k \cup e_k,
k=1,\ldots,m$, where $\widetilde{\Delta}_k$ is a bounded interval
of the real line, $\sigma_k^{\prime} > 0$ a.e. on
$\widetilde{\Delta}_k$, and $e_k$ is a set without accumulation
points in $\mathbb{R} \setminus \widetilde{\Delta}_k$. Let
$\Lambda \subset \mathbb{Z}_{+}^{m}(*)$ be a sequence of distinct
multi-indices such that  $\displaystyle{\max_{{\bf n} \in \Lambda
}(\max_{k =1,\ldots,m}mn_{k} -  |{\bf n}| ) }< \infty$. Let us
assume that there exists $l\in \{1,\ldots,m\}$ and a fixed
permutation $\tau$ of $\{1,\ldots,m\}$ such that for all ${\bf n}
\in \Lambda$ we have that ${\bf n},{\bf n}_l \in
\mathbb{Z}_{+}^{m}(*,\tau)$.  Let $\{q_{{\bf n},k} = \kappa_{{\bf
n},k}Q_{{\bf n},k} \}_{k=1}^m ,{\bf n}\in {\Lambda},$ be the
systems of orthonormal polynomials as defined in
$(\ref{eq:orton})$ and $\{K_{{\bf n},k}\}_{k=1}^m ,{\bf n}\in
{\Lambda},$ the values given by $(\ref{eq:K})$. Then, for each
fixed $k = 1,\ldots,m,$ we have
\begin{equation} \label{eq:xa}\lim_{{\bf n}\in
{\Lambda}}\frac{\kappa_{{\bf n}_l,k}}{\kappa_{{\bf n},k}}=
\kappa^{(l)}_k\,,
\end{equation}
\begin{equation} \label{eq:xk}
\lim_{{\bf n}\in {\Lambda}}\frac{K_{{\bf n}_l,k}}{K_{{\bf n},k}}=
\kappa^{(l)}_1\cdots\kappa^{(l)}_k\,,
\end{equation}
and
\begin{equation} \label{eq:xb}
\lim_{{\bf n}\in {\Lambda}}\frac{q_{{\bf n}_l,k}(z)}{q_{{\bf
n},k}(z)}= \kappa^{(l)}_k \widetilde{F}_k^{(l)}(z), \qquad z \in
\mathcal{K} \subset {\mathbb{C}} \setminus \supp{\sigma_k^{k-1}}
\,,
\end{equation}
where
\begin{equation} \label{eq:xc} \kappa^{(l)}_k =
\frac{c_{k}^{(l)}}{\sqrt{c_{k-1}^{(l)}c_{k+1}^{(l)}}}\,, \qquad
c_{k}^{(l)} = \left\{
\begin{array}{ll}
(F^{(l)}_k)^{\prime}(\infty)\,, & k=1,\ldots,\tau^{-1}(l) \,, \\
F^{(l)}_k(\infty)\,, & k= \tau^{-1}(l)+1,\ldots,m\,,
\end{array}
\right.
\end{equation}
and the $F^{(l)}_k$ are defined by $(\ref{eq31})$.
\end{co}

\begin{proof} By  Theorem  \ref{teofund}, we have limit in
(\ref{eq:var}) along the whole sequence $\Lambda$. Reasoning as in
the deduction of formulas (\ref{eq:otraa}) and (\ref{eq:otra}),
but now in connection with orthonormal polynomials (see Theorems 1
and 2 of \cite{BCL}),   it follows that
\[
\lim_{{\bf n}\in \Lambda} \frac{q_{{\bf n}_l,k}(z)}{q_{{\bf
n},k}(z)}= \left\{
\begin{array}{ll}
(S_k\varphi_k)(z) \,, & k = 1,\ldots,\tau^{-1}(l)\,, \\
S_k(z) \,, & k = \tau^{-1}(l)+1,\ldots,m\,,
\end{array}
\right. \quad \mathcal{K} \subset \mathbb{C}\setminus
\supp{\sigma_k^{k-1}}\,,
\]
where $S_k$ is defined in (\ref{eq:Se}). This formula, divided by
(\ref{eq:otraa}) or (\ref{eq:otra}) according to the value of $k$
gives
\[ \lim_{{\bf n}\in {\Lambda}}\frac{\kappa_{{\bf n}_l,k}}{\kappa_{{\bf n},k}}= \sqrt{\omega_k} =
\frac{c_k}{\sqrt{c_{k-1}c_{k+1}}}\,,
\]
where $\omega_k$ is defined in (\ref{eq:xd}), and the $c_k$ are
the normalizing constants found in Lemma 3.1 solving the linear
system of equations (\ref{eq:bb}) which ensure that
\[ F_k^{(l)} \equiv c_k \widetilde{F}_k^{(l)}\,, \qquad k=1,\ldots,m\,,
\]
with $F_k^{(l)}$ satisfying (\ref{sisec}) and thus given by
(\ref{eq31}). Since $(\widetilde{F}_k^{(l)})^{\prime}(\infty) = 1,
k=1,\ldots,\tau^{-1}(l),$ and $(\widetilde{F}_k^{(l)})(\infty) =
1, k=\tau^{-1}(l)+1,\ldots,m,$ formula (\ref{eq:xa}) immediately
follows with $\kappa_k^{(l)}$ as in (\ref{eq:xc}).

From the definition of $\kappa_{{\bf n},k}\,,$ we have that
\[ K_{{\bf n},k} = \kappa_{{\bf n},1}\cdots\kappa_{{\bf n},k} \,.
\]
Taking the ratio of these constants for the multi-indices ${\bf
n}$ and ${\bf n}_l$ and using (\ref{eq:xa}), we get (\ref{eq:xk}).
Formula (\ref{eq:xb}) is an immediate consequence of (\ref{eq:xa})
and (\ref{eq30}).
\end{proof}

{\bf Acknowledgments.}  G. L\'opez Lagomasino received support
from BFM 2003-06335-C03-02, INTAS 03-516637, and UC3M-MTM-05-033.
A. L\'opez Garc\'ia received support from  UC3M-MTM-05-033.


\begin{thebibliography}{99}
\bibitem{AptLopRoc} A.I. Aptekarev, G. L\'opez Lagomasino, and I.A. Rocha,
{\em Ratio Asymptotic of Hermite-Pade orthogonal polynomials for
Nikishin systems,} Mat. Sb.  {\bf 196} (2005), 1089-1107.
\bibitem{BCL} D. Barrios, B. de la Calle, and G. L\'opez Lagomasino, {\em
Ratio and relative asymptotics of polynomials orthogonal with
respect to varying Denisov type measures,} J. of Approx. Theory
{\bf 139} (2006), 223-256.
\bibitem{BernardoGuillermo1} B. de la Calle Ysern, G. L\'{o}pez Lagomasino, {\em Weak
Convergence of Varying Measures and Hermite-Pad\'{e} Orthogonal
Polynomials,} Constr. Approx. 15 (1999), 553-575.
\bibitem{Denisov} S. A. Denisov, {\em On Rakhmanov's theorem for Jacobi
matrices,} Proc. Amer. Math. Soc. 132 (2004), 847-852.
\bibitem{FidLop} U. Fidalgo and G. L\'{o}pez Lagomasino,  {\em On perfect
Nikishin systems,} Comp. Methods in Function Theory, {\bf 2}
(2002), 415-426.
\bibitem{FidLop2} U. Fidalgo and G. L\'{o}pez Lagomasino, {\em
Rate of convergence of generalized Hermite-Pad\'e approximants of
Nikishin systems,}  Constr. Approx. {\bf 23} (2006), 165-196.
\bibitem{FidLop3} U. Fidalgo and G. L\'{o}pez Lagomasino, {\em
 General Results on the Convergence of Multipoint Hermite-Pad\'e
 Approximants of Nikishin Systems},  Constr. Approx. (accepted).
\bibitem{GonRakhSor} A.A. Gonchar, E.A. Rakhmanov, and V.N. Sorokin,
{\em Hermite-Pad\'{e} for systems of Markov-type functions,} Mat.
Sb. {\bf 188} (1997), 33-58.
\bibitem{Krein} M.G. Krein and A.A. Nudelmann, {\em The Markov
Moment Problem and Extremal Problems,} Transl. of Math.
Monographs, Vol. {\bf 50}, Amer.\ Math.\ Soc., Providence, R.I.,
1977.
\bibitem{kn:Gui1}
G. Lopes [G. L\'opez Lagomasino], {\em On the asymptotic of the
ratio of orthogonal polynomials and convergence of multipoint
Pad\'e approximants,}    Math. USSR Sb.  {\bf 56} (1987), 207-220.
\bibitem{kn:Gui2}
G. Lopes [G. L\'opez Lagomasino], {\em Convergence of Pad\'e
approximants of Stieltjes type meromorphic functions and
comparative asymptotic of orthogonal polynomials,}    Math. USSR
Sb.  {\bf 64} (1989), 207-227.
\bibitem{kn:Nev1}
P. Nevai, {\em Weakly convergent sequences of functions and
orthogonal polynomials,}   J. of Approx. Theory  {\bf 65} (1991),
322-340.
\bibitem{NT}
P. Nevai and V. Totik, {\em Denisov's theorem on recurrence
coefficients,}   J. of Approx. Theory {\bf 127} (2004), 240-245.
\bibitem{kn:Nikishin}
E.M. Nikishin, {\em On simultaneous Pad\'{e} approximations,}
Math. USSR Sb. {\bf 41} (1982), 409-426.

\bibitem{kn:Rak1}
E. A. Rakhmanov, {\em On the asymptotic of the ratio of orthogonal
polynomials,}    Math. USSR Sb.  {\bf 32} (1977), 199-213.

\bibitem{kn:Rak2}
E. A. Rakhmanov, {\em On the asymptotic of the ratio of orthogonal
polynomials II,}     Math. USSR Sb.  {\bf 46} (1983), 105-117.

\bibitem{kn:Rak3}
E. A. Rakhmanov, {\em On asymptotic properties of orthogonal
polynomials on the unit circle with weights not satisfying
Szeg\H{o}'s condition,}    Math. USSR Sb.  {\bf 58} (1987),
149-167.
\end{thebibliography}
\end{document}